%% file: main.tex
\documentclass[a4paper]{article}
\usepackage[utf8]{inputenc}
\usepackage[T1]{fontenc}
\usepackage[french,english]{babel}
\usepackage{soul}
\usepackage[pdftex]{graphicx}
\usepackage{geometry}
\usepackage{tabularx}
\usepackage{float}
\usepackage{chngpage}

\usepackage{tikz}
\usepackage{subfig}
\usetikzlibrary{shapes,arrows,shadows,patterns,decorations.markings,chains}
\usetikzlibrary{positioning}
\usetikzlibrary{calc}

\usepackage[]{algorithm2e}

\usepackage{hyperref}
\usepackage{authblk}

\usepackage{amsfonts}
\usepackage{amsmath}
\usepackage{amsthm}
\usepackage{amsmath}
\usepackage{amssymb}
\usepackage{mathrsfs}
\usepackage{booktabs}
\usepackage{siunitx}
\usepackage{thmtools}
\usepackage{ulem}
\usepackage{mathtools}

\DeclarePairedDelimiter\ceil{\lceil}{\rceil}

\usepackage[explicit]{titlesec}
\usepackage{xcolor}
\usepackage{times}
\usepackage{tikz}
\usepackage{lipsum}
\definecolor{mygreen}{rgb}{0.05, 0.5, 0.06}

\usepackage[title,titletoc]{appendix}

\declaretheoremstyle[
    bodyfont=\normalfont\color{red},
    headfont=\color{red}
]{styleattention}

\declaretheoremstyle[
    spacebelow=1em
]{styleremarque}

\declaretheoremstyle[
    spaceabove=-6pt,
    spacebelow=6pt,
    headfont=\normalfont\bfseries,
    bodyfont=\normalfont,
    postheadspace=1em,
    qed=$\Box$,
]{mystyle}

\declaretheorem[thmbox=M,numberwithin=section,title=Definition]{definition}
\declaretheorem[thmbox=M,sibling=definition]{proposition}
\declaretheorem[thmbox=M,sibling=definition,title=Corollary]{corollaire}
\declaretheorem[thmbox=M,sibling=definition,title=Theorem]{theoreme}
\declaretheorem[thmbox=M,sibling=definition,title=Lemma]{lemme}

\declaretheorem[style=styleremarque,sibling=definition,title=Remark]{remarque}

\declaretheorem[name={}, style=mystyle, unnumbered, title=Proof]{preuve}

\input{command.tex}

\title{Solving Abel integral equations by regularisation in Hilbert scales}
\author{Cecile Della Valle\footnote{Universit\'e de Paris, FP2M, CNRS FR 2036, MAP5 UMR 8145, F-75006 Paris, France. } \hspace{0.05cm}, \hspace{0.5cm}
Camille Pouchol$^\ast$
}
\date{}

\begin{document}
\maketitle
\begin{abstract}
Integral operators of Abel type of order $a>0$ arise naturally in a large spectrum of physical processes. Their inversion requires care since the resulting inverse problem is ill-posed. The purpose of this work is to devise and analyse a family of appropriate Hilbert scales so that the operator is ill-posed of order $a$ in the scale. We provide weak regularity assumptions on the kernel underlying the operator for the above to hold true. Our construction leads to a well-defined regularisation strategy by Tikhonov regularisation in Hilbert scales. We thereby generalise the results of Gorenflo and Yamamoto for $a<1$ to arbitrary $a>0$ and more general kernels.  Thanks to tools from interpolation theory, we also show that the a priori associated to the Hilbert scale formulates in terms of smoothness in usual Sobolev spaces up to boundary conditions,
and that the regularisation term actually amounts to penalising derivatives. 
Finally, following the theoretical construction, we develop a comprehensive numerical approach, where the a priori is encoded in a single parameter rather than in a full operator. 
Several numerical examples are shown, both confirming the theoretical convergence rates and showing the general applicability of the method.
\end{abstract}

\section{Introduction}

\paragraph{Inverse Problem.}
The context of the paper is the inversion of one-dimensional Abel operators of the form
\begin{equation}
\label{def:Abel}
T_a x(t) = \int_0^t (t-s)^{a-1} k(t,s) x(s) \d s 
    \; .
\end{equation}
where $a>0$, and $k$ is a kernel satisfying some appropriate regularity conditions. 
Such operators come up naturally in various physical applications. 
%

For $a<1$, these operators are ubiquitous because they are at the core of fractional dynamical equations. Hence, they play a major role in modelling natural phenomena such as diffusion processes~\cite{mainardi2007}, reaction kinetics of proteins~\cite{glockle1995}, viscoelastic materials~\cite{yang2016}, the physics of surface-volume exchange~\cite{evans2017} to name but a few. They may even be found in applications to psychology~\cite{song2010}.

The specific case where $a=1/2$ is even more broadly studied. In the latter case, the Abel integral $T_a x$ stands for the radial distribution of some
spherically or cylindrically symmetric quantity, such cases arising in modelling plasmas~\cite{merk2013} and flames~\cite{aakesson2008}, in tomography~\cite{Dribinski2002}, or in the so-called \textit{star cluster problem}~\cite{Kosarev1980}. Typically, the inverse problem consists in reconstructing a distribution
of a two-dimensional or three dimensional function from measurements of the projection of these quantities onto a given axis.

The case where $a \geq 1$ can also be found in a variety of applications. 
In hydraulics~\cite{higdon2003inference}, for instance, an Abel integral with $a=3/2$ relates the shape of a notch of a weir and its
flow rate. For $a=2$ or $a=3$, such operators emerge from polymer clustering problems as in~\cite{armiento2016estimation} and for $a=4$ in cristallisation processes~\cite{brivadis2020}.
In these two clustering problems, the experimenter measures the time-evolution of some moment of the polymer distribution. 
Assuming, as is commonly done, that the distribution solves a transport equation with constant or time-varying speed, the inverse problem to be solved belongs to the class~\eqref{def:Abel}.

%
In many of the above applications, 
the functions of interest are smooth functions, such as Gaussian or Gaussian alike. 
Such information of regularity may be taken into account when it comes to improving the inversion strategy in the definition of the prior.

\paragraph{Regularisation strategy.} 
In the present work, we aim at solving this ill-posed problem 
by means of a regularisation strategy of Tikhonov type. More precisely, 
we penalise derivatives of the function we wish to reconstruct. 
Such strategies are commonly used regardless of the operator. At least formally, they are expected to be suitable when the a priori is formulated in terms of smoothness, such as $x$ belonging to some Sobolev space $H^q$, $q>0$.

More specifically, for a measurement $\yd \in L^2(0,1)$ corrupted by noise, we solve the inverse problem by computing
\begin{equation}
    \label{def:AbelInv}
    \xda = \underset{x}{\text{argmin}} \;  
    \| T_a x- \yd \|^2
    + \alpha \|L x\|^2
    \; ,
\end{equation}
where the norm is that of $L^2(0,1)$. Here, and informally at this stage:
\begin{itemize}
\item $\alpha > 0$ is a regularisation parameter, 
\item $\delta> 0$ stands for the noise level,
\item $L$ is a smoothing operator and formally represents the $p$-th derivative of $x$, which requires that $x$ has $p$ derivatives in $L^2$, \textit{i.e.} $x \in H^p$.
\end{itemize}
Regularising by a method of Tikhonov type with underlying smoothing operator $L$ may be studied in the framework of the associated Hilbert scale $(X_p)_{p \in \R}$, with corresponding norms $\|\cdot\|_p$, \textit{i.e.}, $\| L x\| = \|x\|_{p}$, when suitable hypotheses on the operator $L^{1/p}$ hold. Provided that $T_a$ is smoothing in this scale, the convergence of the method is well understood in a very generic framework~\cite{natterer1984, tautenhahn1996}, with infinite smoothing~\cite{mair1994}, and when it comes to finely tuning the regularisation parameter $\alpha$~\cite{Neubauer1988alpha}. 

Although penalising derivatives is common practice, there is no reason that this efficiently achieves the inversion of a given operator $T$. When trying to apply the general framework of Hilbert scales to penalising derivatives and formulating the a priori in terms of smoothness, the following difficulties arise.
\begin{itemize}
\item One needs to build a suitable operator $L$ so that both
\begin{itemize}
\item[(i)] penalising the derivatives of order $p$ is equivalent to penalising the norm $\|\cdot\|_p$,
\item[(ii)] the a priori $x \in X_q$ formulates in terms of usual smoothness assumptions $x \in H^q$.
\end{itemize}
\item Once this is done, the operator $T$ must be shown to be smoothing of some order in the scale~$(X_p)_{p\in \R}$.
\end{itemize}

\paragraph{State of the art.}
Part of this program has been successfully carried out in the works~\cite{Ang1992, Gorenflo1999} in the case where $a \leq 1$. 
The authors show how the Laplace operator associated with appropriate boundary conditions allows one to build 
a well-adapted Hilbert scale~$(X_p)_{p\in \R}$.
However, a series of important questions remains unanswered.

First, the constructed Hilbert scale is suitable only for $a \leq 1$: the operator $T_a$ is not smoothing in the scale constructed in the aforementioned works whenever $a>1$.

Second, the Hilbert scale of~\cite{Gorenflo1999} has only been partially characterised and the link between prior $x \in X_q$ and regularity $x \in H^q$ has only been established for $q \leq 1$, which leaves out any stronger but realistic smoothness assumption, that is when $x$ has more than one derivative in $L^2$.

Third, the efficient numerical implementation of such an approach is up to our knowledge yet to be discussed. At first glance, if one goes from penalising the first derivative to penalising the second, a significant part of the code must be changed. Also, penalising high-order derivatives leads to cumbersome finite difference approximations. Finally, it is not clear how to penalise fractional derivatives.

\paragraph{Main contributions.}
The goal of the present paper is to bridge these gaps. Our contributions may be summed up as follows.

\textit{Construction of an appropriate Hilbert scale.}
We build a family of Hilbert scales, indexed by an integer parameter $\sco$. We show that  regularising in this Hilbert scale exactly amounts to penalising derivatives, and that the a priori $x \in X_q$ corresponds to $x \in H^q(0,1)$ up to some boundary conditions at $t=0$ and $t=1$. 

\textit{Smoothing properties of the operator in the scale.} We then show that the resulting Hilbert scale is suited to the operator $T_a$ when one picks $\sco= \ceil{a}$, assuming enough regularity for the kernel $k$: with this choice, $T_a$ is smoothing of order $a$ in the scale. In fact, we provide two criteria, one which follows and generalises the method of proof of~\cite{Gorenflo1999}, another one of a more functional analytic flavour.  

\textit{Efficient numerical implementation.}
Working around the difficulties mentioned above, we instead closely follow the Hilbert scale at the discrete level. This framework itself advocates for tuning a single parameter $p \geq 0$ standing for which derivative is being penalised. 
Indeed, we explain how a single matrix has to be computed and raised to the chosen power $p$. We confirm our theoretical results and illustrate the flexibility of the approach for various problems involving Abel operators. 



\paragraph{Tools and methods.}

The literature features two main techniques 
when it comes to studying convergence rates for Tikhonov-type regularisation. 

A first category builds upon spectral decompositions and explicit calculations. 
For example, 
one finds results in the case of the Abel integral for $a\leq 1$ in~\cite{groetsch2007}. 
However, 
explicit calculations to build the resulting so-called \textit{filters} are out of reach when $a$ becomes large.

The second family of methods relies on the construction of an adapted Hilbert scale. 
That is the case of the work of~\cite{Gorenflo1999}  for $a \leq 1$. Let us also mention the work~\cite{egger2014}, for very specific cases when $a>1$, where
some simplifications inherent in the problem allow the authors to conclude.

We adopt the latter strategy, but the proof of our main results cannot be carried out as a mere generalisation of~\cite{Gorenflo1999} which heavily relies on the explicit eigensystem of the Laplace operator (with the appropriate boundary conditions). Instead, the operator we need to work with is defined as some possibly higher power of the Laplacian together with suitable boundary conditions. The eigensystem of the resulting operator becomes intractable as $a$ increases, as evidenced by~\cite{bottcher2006} or~\cite{ilin2003}.

Instead, we make extensive work of interpolation theory. Typically, we prove results for specific integer values for which we may directly perform computations such as integration by parts, and then extend the results to fractional values by interpolation. The latter step requires knowledge of interpolation spaces between some standard Sobolev spaces.
For relatively simple cases, the article~\cite{fujiwara1967} provides some results, but the present work requires the more advanced results given in~\cite{guidetti1991}, where general Besov spaces and boundary conditions are treated.

Once the family of Hilbert scales is constructed, we establish that the operator $T_a$ is injective and smoothing of order $a$ in the scale given by to $r = \ceil{a}$, provided that $k$ does not vanish on the diagonal $s=t$ and is sufficiently smooth. We provide two approaches to establish the result which may be complementary depending on the kernel $k$. The first one follows the approach of~\cite{Gorenflo1999} in establishing a suitable factorisation of the operator $T_a$, the second relies on an alternative factorisation together with Peetre's Lemma~\cite{Peetre1961}, but requires $a \notin \N$ and the a priori assumption that $T_a$ is injective.
Then, we may rely on Natterer's Theorem~\cite{natterer1984} to compute the rate of convergence the chosen method has.

\paragraph{Outline of the paper.} 
First, we set up the theoretical framework required for our work in Section~\ref{sec:prem}, \textit{i.e.}, that of fractional Sobolev spaces $H^s(0,1)$, the theory of interpolation of Banach spaces and some results on fractional powers of operators.
Section~\ref{sec:hscale} is devoted to constructing the (integer-indexed) family of  Hilbert scales and identifying it with usual Sobolev spaces.
The next section, Section~\ref{sec:theoremTa}, then provides the main result that the Abel operator is smoothing of order $a$ in the appropriately chosen scale. The convergence of the method is then obtained through a direct application of Natterer's Theorem~\cite{natterer1984}. 
Finally, Section~\ref{sec:numerique} consists of a thorough discussion of how to apply the approach in practice, together with numerical simulations in several contexts involving Abel operators.

%

\section{Mathematical background}
\label{sec:prem}
We introduce the spaces we will be dealing with, namely fractional Hilbert spaces. We also cover the bits of interpolation theory of Hilbert spaces that will be needed throughout. 

We shall always work with spaces of complex-valued functions defined on the interval $(0,1)$. As usual, $\bar{x}$ denotes the complex conjugate of $x\in \C$.
The norm $\|\cdot\|$ and scalar product $(\cdot, \cdot)$ without subscript will refer to the $L^2(0,1)$-norm and scalar product, respectively. The notation $\|\cdot\|$ will also refer to the operator norm of bounded operators from $L^2(0,1)$ onto $L^2(0,1)$.
The identity operator over $L^2(0,1)$ will be referred to as $\1$.

For two Hilbert spaces $\X$ and $\Y$ endowed with respective norms $\|\cdot\|_{X}$ and $\| \cdot\|_\Y$, we will write 
\[
\X \approx \Y
\; ,\]
to indicate that these spaces are topologically equal, \textit{i.e.}, when $\X = \Y$ and the norms $\|\cdot\|_{\X}$ and $\| \cdot\|_\Y$ are equivalent.


\subsection{Fractional Sobolev spaces}
For $k \in \N$, the notation $H^{k}(0,1)$ stands for the usual Sobolev space $W^{k,2}(0,1)$ of functions having $k$ derivatives in $L^2(0,1)$, endowed with the norm 
\[
\|x\|_{H^k(0,1)}^2 := \sum_{j=0}^{k} \|x^{(j)}\|^2 
\; .
\]
For any $\theta \in (0,1)$, the fractional Hilbert space $H^\theta (0,1)$,
is defined as follows
\[
H^\theta (0,1) = 
\left\{
x \in L^2(0,1) \; \text{s.t.} \; \frac{|x(t)-x(s)|}{|t -s|^{1/2 + \theta}} \in L^2((0,1)^2))
\right\} 
\; , 
\]
equipped with the norm
\[
\| x \|_{H^\theta(0,1)}^2 = \|x\|_{L^2(0,1)}^2 + |x|_{H^\theta(0,1)}^2
 \; ,
 \]
 where $| \cdot |_{H^\theta(0,1)}$ is the Gagliardo semi-norm
 \[
 |x|_{H^\theta(0,1)}^2 = 
 \int_0^1 \int_0^1 \frac{|x(t)-x(s)|^2}{|t -s|^{1 + 2\theta}} \d s \d t
 \; .
 \]
 Then, in the case where $s = k + \theta$, $k$ a positive integer and $\theta \in (0,1)$,
the fractional Hilbert space correspond to functions $x$
 whose distributional derivative $x^{(k)}$ belongs to $H^\theta(0,1)$, \textit{i.e.}, 
 \[
H^s (0,1) = 
\left\{
x \in H^k(0,1) \; \text{s.t.} \; x^{(k)} \in H^\theta(0,1)
\right\},
\;
\]
endowed with the norm
\[\|u\|_{H^{s}(0,1)}^2 = \|u\|_{H^k(0,1)}^2 + \|u^{(k)}\|_{H^{\theta}(0,1)}.\]
We insist that the semi-norm (denoted $|\cdot|_{H^s(0,1)}$) of a function $x \in H^s$ refers 
\begin{itemize}
\item to $\|x^{(s)}\|$ when $s$ is an integer,
\item to the the Gagliardo semi-norm $|x|_{H^{\theta}(0,1)}$ when $s$ is not an integer, with $\theta$ denoting its fractional part. 
\end{itemize}
Throughout, whenever the context is clear, we shall drop the reference to the interval $(0,1)$ and use the notation $L^2$, $H^p$ for any $p >0$.

\subsection{The $K$-method for interpolating Hilbert spaces}

Let $X$ and $Y$ be two separable Hilbert spaces with $X$ continuously and densely embedded into~$Y$.
The $K$-interpolation method (which is among the so-called real interpolation methods) is defined as follows: for $t>0$, $y \in Y$, we let
\[K(t,y) := \left(\inf_{x \in X}  \|x\|_X^2 + t^2 \|y-x\|_Y^2\right)^{1/2} \; .\]
For $\theta \in (0,1)$, we let
\[\|y\|_\theta^2 := \int_0^{+\infty} t^{-(2\theta+1)} K(t,y)^2 \,\d t \; .\]
The interpolation spaces are then defined by
\[[X,Y]_\theta := \{ y \in Y, \, \|y\|_\theta < \infty\},\]
endowed with the norm $\|\cdot\|_\theta$.

For our purpose, we will need the following result: interpolating between Hilbert spaces which belong to a Hilbert scale, leads to the expected intermediate space. 
For completeness, we provide a proof of this result.
\begin{lemme}
\label{interpol_scale}
Let $(X_p)_{p \in \R}$ be a Hilbert scale generated by a strictly positive, self-adjoint operator~$D$ with $D^{-1} : X_0 \rightarrow X_0$ compact.
Then for all $0 \leq r < s $, for all $\theta \in (0,1)$, we have \[[X_{s}, X_{r}]_\theta \approx X_{(1-\theta) s + \theta r}.\]
\end{lemme}
Let $(\mu_n , u_n)$ denote an eigensystem for the operator $D$, with $\mu_n>0$, $\mu_n \rightarrow +\infty$ as $n \rightarrow +\infty$ and $(u_n)_{n \in \N}$ an orthonormal basis of $L^2$. Such a decomposition exists by compactness of~$D^{-1}$. Then, recall that the Hilbert scale is characterised for $p \geq 0$ by
\[X_p = \left\{x \in L^2, \,\sum_{n=0}^{+\infty} \mu_n^{2p} \left|(x,u_n)\right|^2 < \infty\right\},\]
and the norm of $x \in X_p$ is given by 
\[\|x\|_p^2 = \sum_{n=0}^{+\infty} \mu_n^{2p} \left|(x,u_n)\right|^2.
\]
\begin{preuve}
Let $\theta \in (0,1)$ be fixed. With the notations above, for $y \in X_{r}$, writing $x = \sum_{n \in \N} x_n u_n$, $y = \sum_{n \in \N} y_n u_n$, the infimum defining $K(t,y)$ may be rewritten as
\begin{align*}K(t,y)^2 & =  \inf_{x \in X_{s}} \|x\|_{s}^2  +  t^2 \|y-x\|_{r}^2 = \inf_{(\mu_n^{s} x_n) \in l^2(\mathbb{N})} \; \sum_{n=0}^{+\infty} \left(  \mu_n^{2s} |x_n|^2  +t^2 \mu_n^{2r} |y_n -x_n|^2 \right).
 \end{align*}
 For each $n \in \N$, the infimum of $x_n \mapsto  \mu_n^{2s} |x_n|^2  +t^2 \mu_n^{2r} |y_n -x_n|^2$ over $\R$ is reached at the value $x_n = \frac{t^2 \mu_n^{2r}} {t^2 \mu_n^{2r} + \mu_n^{2s}}y_n$, and hence equals $ \frac{t^2 \mu_n^{2r} \mu_n^{2s}} {t^2 \mu_n^{2r} + \mu_n^{2s}}|y_n|^2$. For this choice of $x_n$, we indeed have $(\mu_n^{s} x_n) \in l^2(\mathbb{N})$ since 
 \[ \mu_n^{s} |x_n| = \mu_n^{s} \frac{t^2 \mu_n^{2r}} {t^2 \mu_n^{2r} + \mu_n^{2s} }|y_n| \sim 
 \mu_n^{s}  \frac{t^2 \mu_n^{2r}}{\mu_n^{2s}} |y_n| = t^2 \mu_n^{2r-s} |y_n| = o(\mu_n^r |y_n|)
 \]
 Hence we end up with
 \[K(t,y)^2 = t^2  \sum_{n=0}^{+\infty}  \frac{\mu_n^{2r}\mu_n^{2s}} {t^2\mu_n^{2r} + \mu_n^{2s} }|y_n|^2 = t^2  \sum_{n=0}^{+\infty} \mu_n^{2r} \frac{1} {1+ t^2 \mu_n^{2(r-s)} }|y_n|^2 .\]
 Now by Fubini's theorem and the change of variable $u =  \mu_n^{r-s} t$, we may compute
\begin{align*}
\|y\|_\theta^2  = &\; \int_0^{+\infty} t^{-(2\theta+1)} K(t,y)^2 \,\d t\\
= & \; \sum_{n=0}^{+\infty} \mu_n^{2r}  |y_n|^2\int_0^{+\infty} t^{1- 2\theta} \frac{1}{1+t^2 \mu_n^{2(r-s)}}  \,\d t \\
= & \; \left(\int_0^{+\infty}  \frac{u^{1- 2\theta}}{1+u^2}  \,\d u\right)\sum_{n=0}^{+\infty} \mu_n^{2r} \mu_n^{2(s-r)(1- \theta)} |y_n|^2 \\
= & \, C \, \|y\|_{(1-\theta)s +\theta r}^2 \;,
\end{align*}
with the constant $C : = \int_0^{+\infty}  \frac{u^{1- 2\theta}}{1+u^2}  \,\d u = \frac{\pi}{2 \sin(\pi \theta)}$. This ends the proof.
 \end{preuve}
 
\subsection{Fractional powers of operators} 
A final important result on fractional powers of operators is worth mentioning.
Let $A$ with dense domain $\mathcal{D}(A)$ be an accretive operator, \textit{i.e.},
\[
\forall x\in \mathcal{D}(A)\, , \quad 
\text{Re} \left(Ax,x \right) \geq 0
\; .
\]
Recall that $A$ is called m-accretive if $A+\lambda \1$ is furthermore surjective for all $\lambda>0$.
Then the fractional powers of $A$ may be defined, see~\cite{fractionalpowersbook2001}.
We will need the so-called \textit{Heinz-Kato inequality}, which states that if some power $r$ of two m-accretive operators compare, then so do their fractional intermediate powers $\theta r$ for all $\theta \in (0,1)$.
\begin{proposition}[Heinz-Kato Inequality]
\label{prop:Heinz-Kato}
Let $A$ and $B$ be m-accretive operators.
If there exists $r>0$, $C>0$ such that $\mathcal{D}(A^r) \subset \mathcal{D}(B^r)$ and 
   \[  \forall x \in \mathcal{D}(A^r), \quad 
   \| B^r x \|_\X \leq C \|A^r x\|_\X,\]
  then 
  for all $\theta \in (0,1)$, $\mathcal{D}(A^{\theta r}) \subset \mathcal{D}(B^{\theta r})$ and
  there exists $C = C(\theta)>0$ such that
   \[  \forall x \in \mathcal{D}(A^{\theta r}), \quad 
   \| B^{\theta r} x \|_\X \leq C \|A^{\theta r} x\|_\X.\]
\end{proposition}
This result was originally proved in~\cite{heinz1961} for $r=1$. In fact, this result has been extended to the more general case of Banach spaces and sectorial operators having bounded imaginary powers in~\cite{roidos2019heinz} (m-accretive operators are sectorial and have bounded imaginary powers, see~\cite{fractionalpowersbook2001}). The proof is based on Theorem 15.28 in~\cite{Kunstmann}, and straightforwardly extends to arbitrary $r>0$.



\section{Construction of the Hilbert scales}
\label{sec:hscale}

This section deals with constructing and characterising the appropriate Hilbert scales, in which the operator $\Ta$ defined by~\eqref{def:Abel} will be projected. 
More precisely, we let $\sco \in \N^*$ be a fixed integer and build a scale $(X_{\sco,p})_{p \in \R}$ indexed by $\sco$. Recall that $r$ will ultimately be chosen as a function of the exponent $a$ appearing in $T_a$ through $\sco = \ceil{a}$.

\subsection{Defining the scales}
For $\sco \in \N^*$, we define
\begin{equation}
\label{def:B}
\left\{
\begin{array}{rl}
B_\sco = & \; (- \Delta)^\sco \\ 
 \mathcal{D} (B_\sco) = & \; \{ x\in H^{2\sco },  \;
               x^{(k)} (1) =0 \; \text{for} \; 0\leq k <\sco \;,
               x^{(k)} (0) =0 \; \text{for} \; \sco \leq k < 2\sco  \}
               \;. 
\end{array}
\right.    
\end{equation}
Such a definition is motivated by the following link with the inverse problem at hand: for $\sco \in \N^*$, we denote
  \begin{equation}
    \label{def:Sr}
        S_\sco x (t) := \int_0^t (t-s)^{\sco -1}  x(s) \d s \;,
\end{equation}
which is nothing but the integral operator $T_a$ with the constant kernel $k=1$.

Indeed, $S_\sco$ and $B_\sco$ are related as follows.
\begin{lemme}
\label{lem:linkBS}
For all $\sco \in \N^*$, $B_r$ defined by~\eqref{def:B} and $S_r$ defined by~\eqref{def:Sr}, there holds
\[
(r-1)!^2 (B_\sco)^{-1} = S_r^\ast S_r
\; , 
\]
as bounded operators in $L^2$.
\end{lemme}

\begin{preuve}
For $x \in \mathcal{D}(B_\sco)$, we may integrate by parts $\sco$ times to uncover
\[
\begin{split}
    S_\sco B_\sco x (t) &= (-1)^\sco \int_0^t (t-s)^{\sco-1}  x^{(2\sco)}(s) \d s \\
            &= (-1)^\sco \left[ (t-s)^{r-1}  x^{(2\sco-1)}(s)\right]_0^t
             +(-1)^\sco (\sco-1) \int_0^t (t-s)^{\sco-2} x^{(2\sco-1)}(s) \d s \\
            &= (-1)^\sco (\sco-1)!  \int_0^t  x^{(\sco+1)}(s) \d s = (-1)^\sco (\sco-1)! \, x^{(\sco)}(t)
\; .
\end{split}
\]
Indeed, since $x$ belongs to $\mathcal{D} (B_\sco)$, 
the boundary conditions it satisfies are such that 
the integrated terms all vanish, and we obtain $S_\sco^*S_\sco B_\sco x = (\sco-1)!^2 x$. Likewise, we easily check that $S_\sco^* S_\sco$ belongs to $\mathcal{D}(B_\sco)$ for all $x \in L^2$, and by integration by parts we find for $x\in L^2$,
\[
B_\sco S_\sco^*S_\sco x = (\sco-1)!^2 \, x
\; .\]
\end{preuve}

The bounded and symmetric operator $S_\sco^\ast S_\sco$ being self-adjoint, so is $B_r$ as an operator from the range of $S_\sco^\ast S_\sco$ (which the previous lemma shows to be precisely $\mathcal{D}(B_\sco)$) into $L^2$. Furthermore, we also obtain that $B_\sco$ is positive, namely  $(B_\sco x,x) \geq 0$ for all $x \in \mathcal{D}(B_\sco)$. It is even strictly positive since one easily checks that $S_\sco$ is injective by differentiating $\sco$ times the equality $S_\sco x= 0$ (see also next section).

Summing up, $B_\sco$ is a densely-defined, self-adjoint and strictly positive operator. Hence, we may define its real powers, each of which generating a Hilbert scale. In particular, we let 
\[D_\sco :=  B_\sco^{1/2\sco}\] and consider the associated Hilbert scale. 
%
\begin{definition}
 For $r \in \N^*$, we define $(\Xkp)_{p \in \R}$ to be the Hilbert scale induced by the operator $(D_\sco, \mathcal{D}(D_\sco))$, with corresponding norms
 \[
 \| D_\sco^p x \|= : \|x\|_{\sco,p}, \qquad x \in \mathcal{D}(D_\sco^p)
 \; .\]
\end{definition}
Note that the operation of $D_\sco$ should be interpreted as "differentiating once". The justification for this convoluted way of differentiating (and with such boundary conditions) being the above relation  between $S_\sco$ and $B_\sco$. 

In particular, since $S_\sco$ is a Hilbert-Schmidt operator from $L^2$ onto $L^2$, $(D_\sco)^{-1} : L^2 \rightarrow L^2$ is compact.

\subsection{Characterising the scales}
We now characterise the Hilbert scale thus constructed through usual (fractional) Sobolev spaces. We recall the Sobolev embedding $H^s \subset C^0$ for $s \geq \frac{1}{2}$~\cite{fractional2012}. 
Hence for a given $p \geq 0$, the pointwise values $x^{(k)}(0)$ 
and $x^{(k)}(1)$ are well-defined for any $k \in \N$, $k \leq p-\frac{1}{2}$ whenever $x \in H^p$.
\begin{proposition}
\label{prop:expaceXpk}
For $\sco \in \N^*$, $p\geq 0 $, let 
\[\tXkp := \left\{  x \in H^p, \;\;
 \text{for} \; k < p-\frac{1}{2}, 
 \; \; \;
 x^{(k)}(0) = 0     \; \text{if} \;\; 0 \leq k [2r]  < \sco,
 \; \; \;
 x^{(k)}(1) = 0    \; \text{if} \;\; \sco \leq k[2r] <2\sco \; 
 \right\}
 \; ,
\]
equipped with the norm
\[
\|\cdot\|_{\tXkp} = | \cdot|_{H^p}
\; .\]
Then for all  $\sco \in \N^*$ and $p \notin \N + 1/2 $, we have
\[\Xkp \approx \tXkp.\]

\end{proposition}
\begin{preuve}
Our proof goes through the following two steps:
\begin{enumerate}
   \item  we prove the result for the specific values $p = 2\sco m$, $m \in \N$, (with equality of norms and not mere equivalence),
    \item we then generalise the result to any $p \geq 0$, $p \notin \N + \frac{1}{2}$, proving the topological equality of $\Xkp$ and~$\tXkp$.
\end{enumerate}
The idea of the proof is summarised in Figure~\ref{figure_proof}.

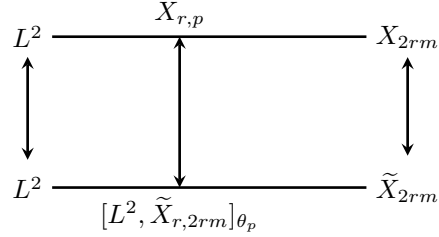
\begin{figure}[h!]
\begin{center}
  \begin{tikzpicture}[
    pizza/.style={black, line width=1, <->, >=stealth},
    anchovy/.style={black, line width=1, <->, >=stealth, shorten >=.1cm},
    peperoni/.style={black, line width=1},
  ]
 \node[] (a1) at (0,0) {$X_{2\sco m}$};
 \node[] (b1) at (0,-2) {$\widetilde{X}_{2\sco m}$};
  \node[anchor=south] (a0) at (-3,0) {$\Xkp$};
 \node[anchor=north] (b0) at (-3,-2) {$[L^2,\widetilde{X}_{\sco ,2 \sco m}]_{\theta_p}$};
 \node[] (a) at (-5,0) {$L^2$};
 \node[] (b) at (-5,-2) {$L^2$};
 \draw[anchovy] (a1) -- (b1);
 \draw[pizza] (a0) -- (b0);
 \draw[anchovy] (a) -- (b);
 \draw[peperoni] (a)--(a1);
  \draw[peperoni] (b)--(b1);
\end{tikzpicture}  
\end{center}
\caption{A schematic idea of the proof.}
\label{figure_proof}
\end{figure}

\paragraph{First step.}
For $p=2 \sco m$, $m \in \N$, 
let us prove by induction on $m$ that the sets $\Xkp$ and $\tXkp$ are equal.
It suffices to prove that $\mathcal{D}(B_\sco^m) = \tXkp$.
For $p = 0$, $\Xkp = \tXkp =L^2$ with the same norm, by definition.

Now assume the equality holds for $p =2 \sco (m-1)$, $m\in \N^*$, and let us address the equality for $p = 2 \sco m$.
\begin{itemize}
 \item Let $x \in \Xkp$, then $x \in \mathcal{D}(B_\sco^m)$ whence $B_\sco x \in \mathcal{D}(B_\sco^{m-1})$.
The induction hypothesis is nothing but $\mathcal{D}(B_\sco^{m-1}) = \tXkpa$, which leads to
$B_\sco x \in H^{p-2r}$ and $x\in H^p$ by elliptic regularity.
Since $\mathcal{D}(B_\sco^{m-1}) = \tXkpa$, we also have
\[
\begin{array}{l}
\text{for} \; 0 \leq r < p-2\sco \; , \; \; 
\begin{cases}
(B_\sco x)^{(k)}(1) = 0 \;& \text{if} \; 0 \leq k \; [2\sco] \leq \sco-1, \\
\;(B_\sco x)^{(k)}(0) = 0 \; &\text{if} \; \sco \leq k\; [2\sco]  \leq 2\sco-1,
\end{cases}
 \\
\\
\Longleftrightarrow 
\text{for} \; 0 \leq k < p-2\sco \; , \; \;
\begin{cases}
x^{(k+2\sco)}(1) = 0 \;& \text{if} \; 0 \leq k \; [2\sco]  \leq \sco-1, \\
\;x^{(k +2\sco)}(0) = 0 \; & \text{if} \; \sco \leq k\; [2\sco]  \leq 2\sco-1,
\end{cases}
 \\
 \\
\Longleftrightarrow 
 \text{for} \; 2\sco \leq k < p \;, \; \; 
 \begin{cases}
  x^{(k)}(1) = 0 \; & \text{if} \; 0 \leq k \; [2\sco] \leq \sco-1, \\
x^{(k)}(0) = 0 \; &\text{if} \; \sco \leq k\; [2\sco] \leq 2\sco-1
 \end{cases}
\; .
\end{array}
\]
Finally, $x\in \mathcal{D}(B_\sco^m)$ does satisfy the required boundary conditions and we indeed have $x~\in~\tXkp$.
\item Conversely, let $x \in \tXkp$, then $B_\sco^k x$ is in $L^2$ for $k \leq m$ and $x \in \mathcal{D}(B_\sco^m)$
and we have $\Xkp =\tXkp$.
\end{itemize}
Now, let us establish the equality of norms. 
For $x \in \Xkp =\tXkp$ with $p = 2 \sco m$, $m \in \N$, we simply write
\[\|x\|_{\Xkp} =  \|D_\sco^{2 \sco m} x\| = \|B_\sco^{m} x\| = \| (- \Delta)^{\sco m}x \| = \|x^{(2 \sco m)}\| = |x|_{H^{2 r m}} =  \|x\|_{\tXkp}.\]

\paragraph{Second step.}
We let $p \geq 0$ be fixed, with $p  \notin \N + \frac{1}{2}$. We pick any $m$ such that $p \leq 2 \sco m$. From the previous step 
\[X_{\sco ,0} = \widetilde{X}_{\sco ,0} = L^2, \qquad  X_{\sco ,2 \sco m} \approx \widetilde{X}_{\sco ,2 \sco m}.\]
Now, we may make use of Theorem 5.1 of~\cite{lions1972}:
for any $\theta \in (0,1)$, the identity $\1$ is continuous as a mapping from the interpolated space
$ [L^2,\Xkp]_\theta $
into the interpolated space
$ [L^2, \tXkp]_\theta $, and conversely. In other words, we have
\[\forall \theta \in (0,1), \quad  [L^2,  X_{\sco ,2 \sco m}]_\theta  \approx [L^2,  \widetilde{X}_{\sco ,2 \sco m}]_\theta.\]
Fixing the value of $\theta$ to $\theta_p :=\frac{p}{2 \sco m}$, we are left with proving that these two spaces are topologically equal to  $\Xkp$ and $\tXkp$, respectively.

In the first case, the result
\[ [L^2,  X_{\sco ,2 \sco m}]_{\theta_p} \approx  X_{\sco, \theta_p 2 \sco m} = \Xkp,\]
is a direct application of Lemma~\ref{interpol_scale}.
On the other hand, identifying $[\widetilde{X}_{\sco,0},\widetilde{X}_{\sco,2 \sco m}]_{\theta_p}$ amounts to interpolating fractional Sobolev spaces with boundary conditions. %
By Theorem 2.7 in~\cite{guidetti1991}, we have for $p \notin \N+1/2 $,
\begin{align*}
    [L^2,\widetilde{X}_{\sco, 2 \sco m}]_{\theta_p} & 
= \widetilde{X}_{\sco, p}.
               \; . 
\end{align*}

\begin{remarque}
Let us mention that in the case $p \notin \N + \frac{1}{2}$, we do not have 
\[ \Xkp \approx \tXkp. \] 
In fact, $\Xkp$ is a (strict) subspace of $\tXkp$~\cite{guidetti1991}.
\end{remarque}
\end{preuve}

\section{Projecting the operator in the Hilbert scale}
\label{sec:theoremTa}

\subsection{Preliminary results}
We now come to our main result, namely that 
$T_a$ is of order $a$ in the scale $\Xkp$, 
where the integer $\sco = \ceil{a} \in \N^*$ denotes throughout this section the smallest integer above $a$.
%
For $a>0$, we extend upon the definition~\eqref{def:Sr} by letting
    \begin{equation}
    \label{def:Sa}
        S_a x (t) = \int_0^t (t-s)^{a -1}  x(s) \d s \;,
    \end{equation}
to which the operator $T_a$ reduces when $k(t,s) = 1$ for all $0 < s < t < 1$.
Then, the integral operator~$S$ is defined for $x\in L^2(0,1)$ through
\[
Sx (t)  = S_1 x(t) = \int_0^t x(\xi) \d \xi
\; .
\]
The family $(S_a)_{a>0}$ satisfies semi-group like properties, see also~\cite{katugampola2011new} for details. Thereafter, $\Gamma$ stands for the usual Euler function \[a \longmapsto \int_0^{+\infty} t^{a-1} e^{-t} \d t.\] 
Let us recall that the fractional powers of $S$ are well-defined, because it is an m-accretive operator. Indeed, since $S$ is bounded, we only need to check that it is accretive to conclude that it is m-accretive. 
For $x \in L^2$, we integrate by parts to find
\[( S x, x) = \int_0^1 x(s)\d s \int_0^1 \bar{x}(s) \d s - (x, S x) \quad \implies \quad \mathrm{Re}( S  x, x) = \frac{1}{2} \left| \int_0^1 x(s)\d s \right|^2  \geq 0. \]
%
\begin{lemme}
\label{lem:semi-group}
For $a \notin \N$, which we write $a =  \sco - 1 + \res$, 
$\sco\in \N^* $ and $0<\res<1$,
we have
 \begin{equation}
\label{prop:factor}
    S_a = \Gamma(\sco)\Gamma(\res) \; S^{\sco-1+\res}  = \Gamma(\sco)\Gamma(\res) \;S^a 
    \;,
\end{equation}
as operators in $L^2$. Moreover, $S_a$ is injective.
\end{lemme}
\begin{preuve}
The case $a \leq 1$,
or equivalently $\sco=1$, 
corresponds to Lemma 5 of~\cite{Gorenflo1999}, 
where it is shown that $S^\res = \Gamma(\res)S_{\res}$.
Hence, the equation~\eqref{prop:factor} holds true 
for $r=1$ and/or $\res =0$
with the convention $S_0 = \1$.
Then, for any $\sco>1$, $x \in L^2$, 
\[
\begin{split}
S_{\sco} x (t) =& \; \int_0^t (t-s)^{\sco-1} x(s)\d s \\
             =& \; \left[(t-s)^{\sco-1} \int_0^s x(\tau) \d \tau \right]_0^t  + \int_0^t (\sco-1) (t-\tau)^{\sco-2}  
             \int_0^s x(\tau) \d \tau \d s\\
             =& \; (\sco-1)S_{\sco-1} S x (t)
             \; .
\end{split}
\]
Then integrating $\sco$ times, we have $S_\sco = \Gamma(\sco) S^{\sco}$.
We have also seen (see~Lemma~\ref{lem:linkBS}) that that $S_\sco$ is injective.
For $a = \sco-1+\res$,
\[
\begin{split}
    S_a x(t) &= \int_0^t (t-s)^{\sco-1+\res-1} x(s) \d s \\
             &= \int_0^t (t-s)^{\sco-1} [(t-s)^{\res-1} x(s)] \d s \\
             &= [(t-s)^{\sco-1} \int_0^t (t-\tau)^{\res-1} x(\tau) \d \tau ]_0^t
                + (\sco-1) \int_0^t (t-s)^{\sco-2} \int_0^t (t-\tau)^{\res-1} x(\tau) \d \tau \d s \\
            &= \; (\sco-1)  S_{\sco-1}S_\res x(t)
            \; .
\end{split}
\]
And therefore
\[ S_{a} = \; (\sco-1) S_{\sco-1} S_\res
= \; \Gamma(\sco) S^{\sco-1} \Gamma(\res) S^\res 
= \; \Gamma (\sco) \Gamma(\res) S^{\sco-1+\res}
\; .\]
Since $S_\res$ is injective (Theorem 4.3 of~\cite{gorenflo1991abel}) and $S_\sco$ is also
injective, $S_a$ is by composition.
\end{preuve}

\subsection{Estimating the order of ill-posedness: first approach}

We now want to establish that $T_a$ is of order $a$ in the Hilbert scale $X_{\sco,p}$ constructed in Section~\ref{sec:hscale}.
The idea is similar to the work~\cite{gorenflo1991abel}: we decompose the operator $T_a$ into a main operator $T_S=k(t,t) S_a$ and a residual operator $T_R$. 
We introduce some useful notations to state our result.
First, we define the open triangle 
\[\Omega := \left\{ (t,s) \in (0,1)^2, \; 0 < s < t < 1 \right\},\]
on which the kernel $k$ is defined. 
%
Secondly, we shall require that $k$ is sufficiently smooth with respect to its second variable and that it does not vanish on the diagonal. More precisely, defining \[g(t,s) :=(k(t,t) - k(t,s)) (t-s)^{a-1}, \quad (t,s) \in \Omega,\]
we assume that
for a.e. $t \in (0,1)$,
\begin{align}
\bullet \quad  & 
  s  \mapsto g(t,s)  \in H^r(0,t) \; , \label{cond:triangle}\\
\bullet \quad    & 
   s \mapsto k(t,s) \;\in  \; \begin{cases}
C^{0,b(a)}(0,t) \; \text{with} \; b(a)>1-\res, \; \; \text{for} \; \; a \leq 1\\
   H^{r-1}(0,t), \; \; \text{for} \; \; a>1
   \end{cases}
  \label{cond:holder}\\
\bullet \quad    & 
       k(t,t) \neq 0 \; . \label{cond:non-nul1}
\end{align}
These hypotheses will be strengthened for our main result to hold true. They are sufficient at this stage to decompose the operator: we can factor $S_a$ out from $T_a$ from the right-hand side.
\begin{lemme}
Under hypotheses~\eqref{cond:triangle}-\eqref{cond:holder}-\eqref{cond:non-nul1}, the decomposition
\begin{equation}
    \label{def:decompRS}
\forall x \in L^2, \quad T_a x (t) =k(t,t) (\1 - R_a) S_a x(t)
\; , 
\end{equation}
holds, where
\begin{equation}
\quad 
R_a x(t) := \int_0^t h(t,s) x(s) \d s
\; 
\end{equation}
Here, denoting $a= \sco-1+\res$, with $\sco = \ceil{a}$, $R_a$ is an integral operator of kernel $h$ vanishing outside of $\Omega$ and given for $(t,s) \in \Omega$ by
\begin{align}
   \bullet \quad & h(t,s) = \frac{1}{k(t,t)} \frac{(-1)^a}{\Gamma(a) }  \frac{\partial^a}{\partial s^a} g(t,s) \; , & a \in \N \; , \label{def:RainN}\\
    \bullet \quad   & h(t,s) =\frac{1}{k(t,t)} \frac{(-1)^\sco }{\Gamma(\sco)} \frac{\sin(\pi \res )}{\pi} 
   \int_s^t (\tau-s)^{-\res} \frac{\partial^\sco}{\partial \tau^\sco} g(t,\tau) \d \tau\; , & a \notin \N \; , \label{def:RaneqN}
\end{align}
with $g$ defined as
\begin{equation}
\label{def:Ra-g}
    g(t,s) :=(k(t,t) - k(t,s)) (t-s)^{a-1}, \quad (t,s) \in \Omega.
\end{equation}
\end{lemme}
\begin{preuve}
For $x \in L^2$, we write
\[
\begin{split}
 T_a x (t) =& \; \int_0^t k(t,s) (t-s)^{a-1} x (s) \d s    \\
       =& \; \underset{  \displaystyle T_S x(t) = \, k(t,t) S_a x (t) }{\underbrace{ k(t,t) \int_0^t(t-s)^{a-1} x (s) \d s}} -
         \underset{\displaystyle T_R x(t)}{ \underbrace{ \int_0^t (k(t,t)-k(t,s) )(t-s)^{a-1} x (s) \d s }} \; ,\\
\end{split}
\]
and by integration by parts, under condition~\eqref{cond:triangle}-\eqref{cond:holder}, all $n \leq \sco $,
\[
\begin{split}
  T_R x(t) =& \; 
            \int_0^t (k(t,t)-k(t,s)) (t-s)^{a-1} x (s) ds \\
        =& \; \left[  (k(t,t)-k(t,s)) (t-s)^{a-1} S x(s) \right]_0^t
        - \int_0^t \frac{\partial}{\partial s} \left(  (k(t,t)-k(t,s)) (t-s)^{a-1} \right) Sx(s) \d s \\
        =& \left[ \frac{\partial^{n-1}}{\partial s^{n-1}} \left(  (k(t,t)-k(t,s))(t-s)^{a-1} \right) S^n x(s) \right]_0^t \\
        & \; \quad + (-1)^{n} \int_0^t \frac{\partial^n}{\partial s^{n}} \left(  (k(t,t)-k(t,s)) (t-s)^{a-1} \right) S^n x(s) \d s \\
       = & \;  (-1)^{n} \int_0^t \frac{\partial^{n} }{\partial s^{n}} \left(  (k(t,t)-k(t,s)) (t-s)^{a-1} \right) S^{n} x(s) \d s \; .
\end{split}
\]
in particular, it is easy to check that the condition~\eqref{cond:holder} is sufficient for the boundary terms to vanish.
In the case of $a \in \N^\star$, or equivalently $a=\sco$, the result is immediate by using $n = \sco$ in the above: formulae~\eqref{def:decompRS} and~\eqref{def:RainN} hold. 

In order to make $R_a$ appear as in formula~\eqref{def:decompRS} when $a= \sco-1+\res$, we still need to factor $S^{\res}$ out from the expression above for $T_R$ with $n = \sco -1$, and justify that the remaining integral operator is well-defined.
The combination of Lemma~\ref{lem:semi-group} and Euler's reflection formula (for $0<\varepsilon<1$, $\Gamma(\varepsilon)\Gamma(1-\varepsilon) = \pi / \sin(\pi \varepsilon)$) leads, for all $x \in L^2$ and $\varepsilon \in (0,1)$, to
\[
S x (t) = \frac{\sin(\pi \varepsilon )}{\pi} \; S_{1-\varepsilon} \; S_\varepsilon x(t) \; .
\]
We introduce the notation $z(t) = S^{\sco-1} x (t)$ and
$f(t,s) = (-1)^{\sco-1} \, \frac{d^{\sco-1}}{ds^{\sco-1}} \left(  (k(t,t)-k(t,s)) (t-s)^{a-1} \right) $, and we have
\[
\begin{split}
  T_R x(t) = & \; \int_0^t f(t,s) z(s) \d s 
         =  \; \int_0^t f(t,s) \frac{\partial}{\partial s}
  \left(  
  \int_{\tau=0}^s z(\tau) \d\tau 
  \right) \d s \\
       = & \;  \frac{\sin(\pi \varepsilon )}{\pi} \int_0^t f(t,s)
          \frac{\partial}{\partial s} \left( 
          \int_0^s (s-\tau)^{-\varepsilon} S_\varepsilon z(\tau) \d \tau 
          \right) \d s \\
    =& \;
    \left[ 
     \frac{\sin(\pi \varepsilon )}{\pi}f(t,s) \int_0^s (s-\tau)^{-\varepsilon} S_\varepsilon z(\tau) \d \tau 
    \right]_0^t\\
    &\; \quad  -  \frac{\sin(\pi \varepsilon )}{\pi}\int_{s=0}^t \frac{\partial}{\partial s} f(t,s) \int_{\tau=0}^s (s-\tau)^{-\varepsilon} 
    S_\varepsilon z(\tau)  \d \tau \, \d s\\
\end{split}
\]
Here we wish to show that the first term cancels out to factorise $T_a$ by $S_a$, which requires $f(t,t)=0$.
For $a \leq 1$ together with the condition~\eqref{cond:holder} we immediately obtain $f(t,t)=0$.
For $a>1$, we may compute explicitly $f(t,t)$ with the Leibniz derivation formula thanks to the condition~\eqref{cond:holder} and check that $f(t,t)=0$ also holds. Hence, we end up with
\[
T_R x(t) 
    = - \frac{\sin(\pi \varepsilon )}{\pi}\frac{1}{ \Gamma(\sco) }
    \int_{s=0}^t \int_{\tau=s}^t (\tau-s)^{-\varepsilon} 
    \frac{\partial}{\partial \tau} f(t,\tau) \;
            \d \tau \, S_{\sco-1+\varepsilon} x (s)   \d s \; .
\]
Let us a posteriori justify the above calculations by proving that the integral 
\[
\int_{s=0}^t \int_{\tau=s}^t (\tau-s)^{-\varepsilon} 
    \frac{\partial}{\partial \tau} f(t,\tau) \;
            \d \tau \, x (s)   \d s 
\; , \]
is well-defined for $x \in L^2$.
Since $ \ell : \tau \mapsto (-\tau)^{-\varepsilon} \delta_{ \tau \leq 0} $ belongs to $L^1(0,t)$ and 
since $k$ satisfies condition~\eqref{cond:triangle},
$\partial_\tau f(t,\tau)$ is in $L^2(0,t)$.
By Young's convolution inequality,
the convolution 
\[
s \mapsto \int_{0}^t (\tau-s)^{-\varepsilon} \delta_{s-\tau \leq 0}  \frac{\partial}{\partial \tau} f(t,\tau)  \, \d \tau 
=  \left( \ell * \frac{\partial}{\partial \tau} f(t,\cdot) \right) (s)
\; ,
\]
belongs to $L^2(0,t)$ .
The integral is therefore the scalar product of two functions of $L^2(0,t)$.

To conclude, we set $\varepsilon = \res$, and we find \eqref{def:decompRS} and \eqref{def:RaneqN} with $T_R = (k(t,t))^{-1} R_a S_a $.
\end{preuve}
Thanks to the above decomposition~\eqref{def:decompRS}, we may exhibit a sufficient condition to compare $T_a$ and~$S_a$. Thus, the two operators are of the same order in the appropriate Hilbert scale.
%
\begin{theoreme}
\label{th:Ta}
Let $a>0$.
Assume that $k$ satisfies~\eqref{cond:triangle}-\eqref{cond:holder}-\eqref{cond:non-nul1}. 
If $k$ is furthermore such that
\begin{align}
\bullet \quad    & 
       c_k^{-1} \leq |k(t,t)| \leq c_k \quad \text{for} \; c_k>0 \; , \label{cond:bound} \\
\bullet \quad  & 
\1-R_a : L^2 \to L^2\; \text{is bounded, invertible}
   \; , \label{cond:Ra-bib}
\end{align}
then $T_a : L^2 \to L^2$ is injective and there exists a constant $c =c(a)$ such that
\[\forall x \in L^2, \quad
c^{-1} \|x\|_{\sco,-a}
\leq
\| T_a x\|_{L^2}
\leq
c \|  x\|_{\sco,-a}
\; .
\]
\end{theoreme}
\begin{remarque}
At first glance, Condition~\eqref{cond:Ra-bib} might seem rather abstract. Since there are numerous sufficient criteria to decide whether such a result holds true (especially coming from Fredholm Theory), we prefer to give the result with such generality and only then to give some workable sufficient conditions ensuring Condition~\eqref{cond:Ra-bib}, see Corollary~\ref{cor:Raexplicit}.
\end{remarque}
\begin{preuve}
\textit{First step: $k=1$}.
In this case, $T_a = S_a$, in which case injectivity has already been established. 
Let us start with the integer case $a =\sco$.
Recalling $(\sco-1)!^2 (B_\sco)^{-1} = S_\sco^\ast S_\sco$ from Lemma~\ref{lem:linkBS}, we have for $x\in L^2$
\[
\begin{split}
    \|S_\sco x\|^2 = (S_\sco^*S_\sco x, x) = (r-1)!^2 (B_\sco^{-1} x, x) =  (r-1)!^2 \|x\|_{\sco,-\sco}^2
    \; ,
\end{split}
\]
Hence, Theorem~\ref{th:Ta} holds true with the constant $c = c(\sco)= (\sco-1)!$

Let $a>0$ be fixed. Letting $\sco = \ceil{a}$, we use the above according to which we have
\[\forall x \in L^2, \quad (r-1)!\| D_\sco^{-\sco} x\|= \| S_\sco x\| \;  \implies \; \| D_\sco^{-\sco} x\|= \| S^\sco x\| \]
using $S_\sco = (r-1)!S^\sco$.
We now aim at applying the Heinz-Kato inequality~\ref{prop:Heinz-Kato}.
Since the operator $D_\sco^{-1}$ is positive and self-adjoint, it is m-accretive. We also know that $S$ is m-accretive. Hence, we may use the Heinz-Kato inequality with $n=\sco$ and $\theta = a/\sco \leq 1$ in Proposition~\eqref{prop:Heinz-Kato}. 
The result is proved because Lemma~\ref{lem:semi-group} shows $S_1^{\theta}S_{r-1}$ and $S_a$ differ only by a multiplicative constant.

\textit{Second step: general $k$}.
For any $k$ under Conditions~\eqref{cond:triangle}-\eqref{cond:non-nul1},
we recall that 
\[
T_a = k(t,t)(\1 - R_a)S_a 
\; . 
\]
Under Conditions~\eqref{cond:triangle}-\eqref{cond:holder}-\eqref{cond:non-nul1}, for $x \in L^2(0,1)$, $R_a x$ is well defined. Moreover, under Condition~\eqref{cond:Ra-bib}, $R_ax$ belongs to $L^2$, and the injectivity of both $\1-R_a$ and $S_a$ yield that of $T_a$.
The hypotheses allow us to bound as follows for $x \in L^2$:
\[
\|T_a x \| \leq c_k \| \1 -R_a \|  \| S_a  x\|\, 
\;,
\]
as well as 
\[
 \| S_a  x\| \leq \, c_k \| (\1 -R_a )^{-1} T_a x \| \leq \, c_k  \| (\1 -R_a )^{-1}  \| \, \| T_a x \|
\;,
\]
and the first step concludes the proof.
\end{preuve}
%
\begin{remarque}
In the proof of Theorem~\ref{th:Ta}, we used the Heinz-Kato inequality~\ref{prop:Heinz-Kato} 
to compare $S_a$ and $D_\sco^{-a}$ with $\sco=\ceil{a}$.
We could in fact have compared $S_a$ with $D_\sco^{-a}$ for any $\sco \geq a$.
Hence, we have actually proved the following result under the hypotheses of the previous theorem: 
\[\text{$T_a$ is of order $a$ in any Hilbert scale $(\Xkp)$ with $\sco \geq a$}.\] It is of course of little interest to take any $r$ larger than $\lceil a \rceil$ in practice if $a$ is known. On the contrary, if there is uncertainty on the value of $a$, say in the form of a weak information $a \in [a_{\min}, a_{\max}]$, then one could (and should) take $r = \lceil a_{\max} \rceil$.
 \end{remarque}

We propose to give sufficient conditions for which 
the condition~\eqref{cond:Ra-bib} on the integral operator $R_a$ is verified, 
which might be more handy to check depending on the inverse problem at hand.
\begin{corollaire}
\label{cor:Raexplicit}
Let $a>0$.
Assume that $k$ satisfies~\eqref{cond:triangle}-\eqref{cond:holder}-\eqref{cond:non-nul1}-\eqref{cond:bound}-\eqref{cond:Ra-bib}. 
If the integral operator $R_a$ of kernel $h$ defined by~\eqref{def:RainN}-\eqref{def:RaneqN} 
satisfies one of the following conditions :
\begin{align}
\bullet \quad    & 
      \|R_a\| < 1 ,  \label{cond:Rainf1}\\
\bullet \quad  & 
\exists \, \gamma \in L^2(0,1) \; \text{ s.t.} \; \; \forall x \in  L^2 \, , \; | R_a x (t)| \leq \int_0^t \gamma (s) |x (s) | \d s   \; ,  \label{cond:intgamma}\\
\bullet \quad  & 
  h  \in L^\infty (\Omega)  \; , \label{cond:hLinfty}
\end{align}
then the conditions of Theorem~\eqref{th:Ta} are satisfied.
\end{corollaire}
\begin{remarque}
Condition~\eqref{cond:hLinfty} certainly implies condition~\eqref{cond:intgamma} (by taking the constant function $\gamma := \|h\|_{L^\infty(\Omega)}$).
We choose to stress condition~\eqref{cond:hLinfty} independently since it may be checked more directly. 

We also note that if $h \in L^2(\Omega)$, $R_a$ is a Hilbert-Schmidt operator, and hence is compact from $L^2$ onto $L^2$. Hence, by the Fredholm alternative, showing that $\1 - R_a$ is bijective is equivalent to showing that it is either injective or surjective.
\end{remarque}
\begin{preuve}
We prove that all conditions imply that $\1 - R_a$ is bounded, invertible with bounded inverse.

\paragraph{Under Condition~\eqref{cond:Rainf1} ,}
since $\|R_a\| <1$, it is standard that $\1-R_a$ is invertible, bounded (with inverse given that then the Neumann series $\sum R_a^n$ of $R_a$).
Consequently, the condition~\eqref{cond:Ra-bib} of Theorem~\ref{th:Ta} holds true.

\paragraph{Under Condition~\eqref{cond:intgamma} , }
the function $R_a x$ is in $L^2$. 
Our reasoning follows the proof of Lemma~1 of~\cite{Ang1992}, which for completeness we repeat below.
Let us prove by induction that 
\[
|R_a^n x(t) | \leq \frac{1}{n!^{1/2} (n-1)!^{1/2}}  \left(\int_0^t \gamma (s)^2 \d s \right)^{n/2} \,
 \left(\int_0^t (t-s)^{n-1} |x(s)|^2  \d s\right)^{1/2} \, ,
\quad n\in \N
\; .
\]
This inequality holds true for $n=1$.
Now for any $n \in \N$, by Hölder's inequality,
\[
\begin{split}
    | R_a^{n+1} x  (t) | \leq & \;  \int_0^t \gamma(s) | R_a^{n} x  (s)| \d s \\
                       \leq & \;  \frac{1}{n!^{1/2}(n-1)!^{1/2}}  \int_0^t \gamma(s)  \left(\int_0^s \gamma (\tau)^2 \d \tau \right)^{n/2} \,
 \left(\int_0^s (s-\tau)^{n-1} |x(\tau)|^2  \d \tau \right)^{1/2} \d s \\
                       \leq & \;  \frac{1}{n!^{1/2}(n-1)!^{1/2}} 
                       \left( \int_0^t \gamma(s)^2  \left(\int_0^s \gamma (\tau)^2 \d \tau \right)^{n} \d s \right)^{1/2} \, \left( \int_0^t \int_0^s (s-\tau)^{n-1} |x(\tau)|^2 \d \tau  \d s\right)^{1/2} \\
                        \leq & \;  \frac{1}{(n+1)!^{1/2}(n)!^{1/2}} 
                       \left( \int_0^t  \left(\int_0^s \gamma (\tau)^2 \d \tau \right)^{n+1} \right)^{1/2} \, \left( \int_0^t (t-s)^{n} |x(s)|^2  \d s\right)^{1/2} \\
                       \leq & \;  \frac{1}{(n+1)!^{1/2}(n)!^{1/2}} 
                        \left(\int_0^t \gamma (s)^2 \d s\right)^{(n+1)/2} \, \left( \int_0^t (t-s)^{n} |x(s)|^2  \d s\right)^{1/2} \;.
\end{split}
\]
Now, if we bound the kernel integral operator $S_n$,
\[
\left( \int_0^1 \int_0^t (t-s)^{n} |x(s)|^2  \d s)^2 \d t \right)^{1/2}
\leq \frac{\| x\| }{\sqrt{n}}
\]
we obtain for all $n \in \N$,
\[ 
\|R_a^n\| \leq \frac{\| \gamma \|^n}{n!} 
\, . \]
Therefore the Neumann series of $R_a$ converges in the operator norm, and $\1 - R_a$ is then invertible with bounded inverse, and the condition~\eqref{cond:Ra-bib} is met.

\end{preuve}
Let us now make condition~\eqref{cond:Rainf1} a bit more explicit.
Since $R_a$ is an integral operator, we can always control its norm by the Hilbert-Schmidt norm (which may or may not be finite) through
\[\|R_a\| \leq \|h\|_{L^2((0,1)^2)} = \|h\|_{L^2(\Omega)}.\]
Hence, a sufficient condition for condition~\eqref{cond:Rainf1} to be satisfied is explicitly given by
\begin{itemize}
    \item for $a \in \N$,
    \[ \frac{1}{\Gamma(a)}\int_\Omega \bigg( \frac{1}{k(t,t)} \frac{\partial^a}{\partial s^a} \left((k(t,t)-k(t,s))(t-s)^{a-1}\right)\bigg)^2 \d s \d t \ <  1
     \; ,
     \]
     \item for $a=\sco-1+\res$, $\res<1$,
     \[ \frac{|\sin(\pi \res)|}{\pi\Gamma(\sco)}\int_\Omega \bigg(
    \frac{1}{k(t,t)}\int_s^t (t-\tau)^{-\res}\frac{\partial^\sco}{\partial \tau^\sco} \left((k(t,t)-k(t,\tau))(t-\tau)^{a-1}\bigg)  \d \tau\right)^2 \d s \d t 
     < 1 
     \; .
     \]
\end{itemize}

\subsection{Estimating the order of ill-posedness: second approach}

We now propose an alternative way to prove the result with less restrictive assumptions, but for the case $a \notin \N$ and assuming that injectivity of $T_a$ has been established independently, and in sufficiently weak spaces. We refer to~\cite{gorenflo1991abel} for some sufficient conditions regarding injectivity for Abel operators in classical $L^p$ spaces, and to Remark~\ref{rem:injective} when it comes to passing from classical injectivity to weaker injectivity (at least for $a<1$). 

Indeed, in the case where $a$ is not integer, $S_{a+\poc}$ can be compared to $S_a$ in the same Hilbert scale provided that $\gamma$ is taken small enough.
This is not the case for $a \in \N$, 
and the underlying reason is that the operators $S$ and $S^*$ do not commute.
Under assumptions~\eqref{cond:triangle}-\eqref{cond:holder}-\eqref{cond:non-nul1}, we go back to the formula~\eqref{def:RaneqN} but for $\varepsilon = \res+\poc$.
We obtain for $x \in L^2$,
\begin{equation}
\label{def:decompRSeta}
    T_a x(t) = k(t,t)S_a x(t) + k(t,t) R_{a,\poc} S_{a+\poc} x(t)\; ,
\end{equation}
with
\begin{equation}
\label{def:Raeta}
     R_{a,\poc} x(t) = c(a,\poc)  \,
    \int_0^t \int_s^t (\tau-s)^{-\res -\poc} \frac{\partial^\sco}{\partial \tau^\sco} g(t,\tau) \d \tau \, x(s) \d s \; ,
\end{equation}
with
$g$ defined by~\eqref{def:Ra-g}
and
\[
c(a,\poc) = \frac{(-1)^\sco}{\Gamma(\sco)} \frac{\sin(\pi (\res + \poc) )}{\pi}  \; .
\]
Here, in contrast to~\cite{Gorenflo1999}, our approach completely diverges. 
Our idea is to view the residual $R_{a,\poc}S_{a+\poc}$ as a compact perturbation of the main part, and to use a corollary of the so-called Peetre's lemma~\cite{Peetre1961}. 
\begin{theoreme}
\label{th:Taeps}
Let $a>0$, $a \notin \N$.
Assume that $k$ satisfies~\eqref{cond:triangle}-\eqref{cond:holder}-\eqref{cond:non-nul1}.
If $T_a$ is injective from $X_{\sco, -a}$ onto $L^2$, and if, for $\gamma$ small enough, 
$R_{a,\poc}$ defined by~\eqref{def:decompRSeta}-\eqref{def:Raeta} satisfies
\begin{equation}
 \label{cond:Retafini}
\|R_{a,\poc}\| < +\infty \; ,   
\end{equation}
then, there exists a constant $c =c(a)$ such that
\[\forall x \in L^2, \quad
c^{-1} \|x\|_{\sco,-a}
\leq
\| T_a x\|_{L^2}
\leq
c \|  x\|_{\sco,-a}
\; .
\]
\end{theoreme}
We emphasise that Theorem~\ref{th:Taeps} requires less demanding conditions, since they do not rely on estimating the norm of the operator $R_{a, \poc}$ (but only that it be finite).
However, by relaxing the assumption~\eqref{cond:Rainf1}, 
another condition naturally emerges, which is that the injectivity of~$T_a$ has to be checked independently since Condition~\eqref{cond:Retafini} is not sufficient to guarantee it.

Before coming to the proof of Theorem~\ref{th:Taeps}, we develop the necessary (functional analytic) tools. In what follows, $C>0$ is a constant which might change from line to line. 
\begin{lemme}
\label{lem:peetre-corollaire}
Given three Hilbert spaces  $X$, $Y$, $Z$ 
such that the injection of~$Y$ into~$Z$ is compact,
and given two bounded operators $S$ et $R$ from~$Y$ to~$X$, 
we define the operator $T := S + R$ which we assume to be injective.
Moreover, if we have
\begin{align}
    & \forall x \in Y, \qquad C^{-1} \|x\|_{Y} \leq \|S x\|_{X} \leq C \|x\|_{Y} \label{hyp1}\\
    & \forall x \in Z, \qquad \|R x\|_{X} \leq C \|x\|_{Z} \label{hyp2}
\end{align}
Then, 
\[
\forall x \in Y, 
\qquad C^{-1} \|x\|_{Y} \leq \|T x\|_{X} \leq C \|x\|_{Y}
\; .
\]
\end{lemme}
\begin{preuve}
The upper bound directly follows from the continuous injection of $Y$ into $Z$,
\[
\forall x \in Y, 
\qquad 
\|T x\|_{X} 
\leq  \|S x\|_{X}+ \|R x\|_{X} 
\leq C \|x\|_{Y} +C \|x\|_{Z} 
\leq C \|x\|_{Y}
\; .
\]
The lower bound, on the other hand, is a consequence of Peetre's lemma, which we recall below.
\begin{lemme}[\cite{Peetre1961}-Lemma 3]
\label{lem:peetre}
Let $X$ and $Y$ be two Banach spaces, and two operators, $T$ bounded injective from $Y$ into $X$,
and $R$ compact from $Y$ into $X$.
We assume that 
\[
\forall x \in Y, 
\qquad 
\|x\|_{Y} \leq C \left(\|Tx\|_{X} + \|Rx\|_{X}\right)
\; .
\]
Then, 
\[
\forall x \in Y, 
\qquad 
\|x\|_{Y} \leq C \|T x\|_{X}
\; .
\]
\end{lemme}
We may write
\[
\forall x \in Y, 
\qquad 
\|S x\|_{X} = \|T x - Rx\|_{X} 
\leq \|T x\|_{X} + \|R x\|_{X}
\; ,
\]
which leads via hypothesis~\eqref{hyp1} to 
\[
\forall x \in Y, 
\qquad 
\|x\|_{Y} 
\leq C\left( \|T x\|_{X} +\|R x\|_{X}\right)
\; .
\]
Now let $(x_n)_{n\in \N}$ be a bounded sequence of~$Y$.
Then, by the compact injection of $Y$ into $Z$, 
$(x_n)_{n\in \N}$ converges in $Z$ upon extraction. 
By the inequality~\eqref{hyp2}, $(Rx_n)_{n\in \N}$ converges in $X$ along the subsequence, which shows that $R$ is compact from $Y$ into $X$. Hence, Lemma~\ref{lem:peetre} concludes the proof.
\end{preuve}
\vspace{0.2cm}

\begin{preuve}
As mentioned in the sketch of proof, we pick $\poc>0$ small enough so that $\ceil{a+\poc} = \ceil{a} = \sco$,
which is possible since $a \notin \N$.
Then, we have
\[
T_a x(t) = k(t,t) S_a x(t)+ R_{a,\poc} S_{a+\poc} x(t)
\; ,
\]
a decomposition which can be shown to hold as in Lemma~\ref{def:decompRS} thanks to~\eqref{cond:triangle}-\eqref{cond:holder}-\eqref{cond:non-nul1}.

The operator $T_a$ is now the sum of two operators, 
a first one whose projection in a Hilbert scale is known thanks to Theorem~\ref{th:Ta}
for the simple case of $k(t,s)=k(t,t)$, 
and a second one that is a compact perturbation of the first one, as we shall see.   

We define $X=L^2$, $Y= X_{\sco,-a}$ and $Z =X_{\sco,-(a+\poc)}$,
the injection of $Y$ into $Z$ is compact since the two spaces belong to the same Hilbert scale.
The operator $S_a$ and $R_{a,\poc}S_{a+\poc}$ are bounded operators from $X = L^2$ to $X = L^2$, respectively thanks to Lemma~\ref{lem:semi-group} and thanks to the condition~\eqref{cond:Retafini}. Hence, they are also bounded operators from $Y = X_{\sco,-a}$ to $X = L^2$ (once uniquely extended) since $L^2 = X_{\sco, 0}$ is densely embedded into $X_{\sco,-a}$.
By assumption, $T_a$ is injective from $Y =X_{\sco,-a}$ onto $X=L^2$.
%
%
By~\eqref{cond:bound}, 
$c_k^{-1} \leq |k(t,t)| \leq c_k $ and using Theorem~\ref{th:Ta}, the first term satisfies
\[
c^{-1} \|x\|_{\sco, -a} \leq \|k(t,t) S_a x \|_{L^2} \leq c \|x\|_{\sco, -a}
\; ,
\]
and the condition~\eqref{hyp1} is met. 

We now evaluate the operator $R=R_{a,\poc}S_{a+\poc}$ in the Hilbert scale $(X_{\sco,p})_{p \in\R}$.
Then, by Theorem~\ref{th:Ta}, there exists $C>0$ such that
\[
C^{-1} \|x\|_{\sco, -(a+\poc)} \leq \| S_{a+\poc} x \| \leq C \|x\|_{\sco, -(a+\poc)}
\; .
\]
Moreover, the operator $R_{a,\poc}$ is assumed to be bounded, 
and we denote $\|R_{a,\poc} \| = M$.
Hence, we have the announced upper bound
\[
\|R_{a,\poc} S_{a+\poc} x \|_{L^2} \leq M \|S_{a+\poc}x\|_{L^2} \leq cM \|x\|_{\sco, -(a+\poc)}
\; .
\]
We have then shown that condition~\eqref{hyp2} holds true.
We may now apply Lemma~\ref{lem:peetre-corollaire}.
and conclude that
for $a \notin \N$, for all $x \in L^2$,
\[
C^{-1} \| x\|_{\sco,a} \leq \| T_a x\|_{L^2} \leq C \| x\|_{\sco,a}
\; .
\]
%
%
\begin{remarque}
\label{rem:injective}
The requirement that $T_a$ be injective from the weak space $X_{\sco,-a}$ into $L^2$ is rather abstract. Since the injectivity of $T_a$ in $L^p$ spaces has thoroughly been studied in the literature, let us give some sufficient conditions such that 
\[T_a :L^2 \to L^2 \text{ injective} \quad \implies \quad T_a :X_{\sco,-a} \to L^2 \text{ injective}.\] 
We work in the setting where $0<a<1$: we leave open the problem of finding comparably simple conditions when $a>1$.

Let $x \in X_{\sco,-a}$ such that $T_a x = 0$. We prove that this equality enforces $x \in L^2$. 
Since decomposition~\eqref{def:decompRS} holds, we have
\[
S_a x(t) = -R_a S_a x(t) \; .
\]
Now, assuming for the moment the following properties
\begin{itemize}
    \item $y \in L^2 \implies R_a y\in H^1$,
     \item $S_a x \in H^a \implies x \in L^2$,
\end{itemize}
we apply the first one to $y =S_a x$, obtaining $y = S_a x \in H^1$. The second property then yields~$x\in~L^2$.
Let us finally discuss the two above properties. The second one is a very general property of $S_a$ that requires no further conditions, and is established in~Theorem 2.1 of~\cite{gorenfloyamamoto}.
The first implication is obtained using the conditions
\begin{align}
\bullet \quad    & 
       \forall s \in (0,1), \; t \mapsto h(t,s) \in H^1(s,1) \; , \label{cond:definition} \\
\bullet \quad  & 
(t,s) \mapsto \frac{\partial}{\partial t} h(t,s) \in L^2(\Omega)
   \; . \label{cond:domination}
\end{align}
Indeed, when these hold, we may use that $h(t,t) = 0$ and differentiate once (in the weak sense) $y = R_a y$ to obtain
\begin{equation}
y'(t)
=
 \int_0^t \frac{\partial}{\partial t} k(t,t) h(t,s) y(s) \mathrm{d} s
\; .
\end{equation}
%
\end{remarque}

%
\end{preuve}
In Appendix~\ref{annexe:proof}, we also elaborate on the particular case where $k$ is analytic with respect to its second variable, which leads to sufficient conditions that may be easier to check in some specific cases. 
\subsection{Tikhonov regularisation}

Let us now return to the solution of the inverse problem associated to the Abel integral.  
We wish to reconstruct $x$ such that 
\[
T_a x = y \; , 
\]
Instead of having access to the exact data $y$, 
we must reconstruct the signal from noisy data $\yd$ 
such that the measurement error is bounded in the $L^2$ norm,
\begin{equation}
\label{def:error}
  \| \yd-y \| \leq \delta   
  \; ,
\end{equation}
We also  are given some a priori regularity about the unknown $x$, which we assume writes
\begin{equation}
\label{def:regx}
    x \in X_{\sco,q}, \quad \|x \|_{\sco,q}\leq M
    \; .
\end{equation}
The Tikhonov-type regularisation method for recovering $x$
consists in solving the minimisation problem
\begin{equation}
\label{def:Lu}
   \min_{x \in X_{r,p}} \mathscr{J} (u)  \qquad \mathscr{J} (u) := \| T_a u - \yd \|_{L^2}^2 + \alpha \| u \|_{\sco,p}^2
\end{equation}
We then have the following convergence theorem, which directly follows from Natterer's theorem~\cite{natterer1984}.
\begin{corollaire}
\label{cor:order}
Under the hypotheses of Theorem~\ref{th:Ta} or Theorem~\ref{th:Taeps}, assume~\eqref{def:error} and~\eqref{def:regx}, take
\begin{equation}
    \label{cond:reg}
    p \geq \frac{q-a}{2} \; ,
\end{equation}
as well as 
\begin{equation}
    \label{choice:reg}
    \alpha = \alpha(\delta) = C \delta^{\frac{2(a+p)}{a+q}} \; ,
\end{equation}
for some $C>0$.
Then, the solution $x^{\delta, \alpha(\delta), p}$ of the minimisation problem~\eqref{def:Lu} with regularisation parameter 
$\alpha(\delta)$ satisfies
\begin{equation}
    \| x^{\delta, \alpha(\delta), p} - x \| \leq c \; \delta^\frac{q}{q+a} M^\frac{a}{a+q} \; 
\end{equation}
for some constant $c>0$.
\end{corollaire}
We repeat here the implication drawn by Natterer in~\cite{natterer1984}: \textit{there is nothing wrong with high order regularisation, even well above ther order of smoothness of the exact solution. The only mistake one can make is to regularise with an order which is too low}.
\begin{remarque}
\label{rem:commute}
Let us also mention the following slight improvement, still from~\cite{natterer1984}:
if one further assumes that
\begin{equation}
    \label{cond:commute}
(T_a^*T_a)^{1/2} \text{ and } D_\sco \text{ commute},
\end{equation}
then Corollary~\ref{cor:order} holds true 
and condition~\eqref{cond:reg} becomes
\[
p \geq \frac{q}{2} -a 
\;.
\]
From the results of Theorem~\ref{th:Ta}, the above hypothesis~\eqref{cond:commute} happens to be satisfied when the kernel $k$ is identically $1$ and $a = \sco$ is a positive integer, \textit{i.e.}, with our notations when $T_a = S_a = S_\sco$. Indeed, recall the equality $B_r^{-1} = S_r^\ast S_r$, which directly entails that $T_a^\ast T_a = S_\sco^\ast S_\sco = B_\sco^{-1}$ commutes with $D_\sco = B_\sco^{\frac{1}{2r}}$.
\end{remarque}
%
\section{Numerical experiments}
\label{sec:numerique}

In this section,
we discuss the numerical and practical solution of an inverse problem related to 
an operator of the form~\eqref{def:Abel} 
by penalising derivatives. 
\subsection{Preliminary remarks}

Corollary~\eqref{cor:order} shows that the quality of the inversion, 
or equivalently of the reconstruction, 
is improved if the function to be reconstructed is smooth and satisfies some boundary conditions. 

In particular, if the unknown is compactly supported inside $(0,1)$, our results simply mean the following: the smoother (in the usual $L^2$ Sobolev sense) the unknown, the better the reconstruction method works.
If however the function is not compactly supported, the order of convergence is controlled by the boundary conditions, even for arbitrarily smooth functions.

Proposition~\ref{prop:expaceXpk} also shows that the minimisation problem 
\begin{equation}
   \min_{x \in X_{r,p}} \| T_a u - \yd \|_{L^2}^2 + \alpha \| u \|_{\sco,p}^2
\end{equation}
is equivalent to 
\begin{equation}
   \min_{x \in X_{r,p}} \| T_a u - \yd \|_{L^2}^2 + \alpha |u|_{H^p}^2,
\end{equation}
upon changing the parameter $\alpha$, and at least for $p \notin \N + \frac{1}{2}$. For $p \in \N$, in particular, this is nothing but penalising the $p$th derivative of $u$ through $\|u^{(p)}\|^2$. 

At this stage, in order to solve the above, we need to  elaborate on how to discretise the $p$th derivative as well as how to deal with the boundary conditions.
As mentioned in the introduction, 
a first approach would be to use finite differences to approximate the derivative~$u^{(p)}$.
Then, for each chosen level or regularisation $p$,
this method would lead to cumbersome computations,
at least when $p \in \N$. 
This has three major drawbacks: 
\begin{itemize}
    \item the code must significantly be changed for each instance of $p$ and may become heavy for large values (recall that we should not refrain from taking $p$ large), 
    \item this does not carry over to the case of fractional $p$,
    \item the boundary conditions are not properly taken into account.
\end{itemize}
In fact, dealing with all three caveats is achieved by closely following the initial formulation 
with the Hilbert scale and underlying operator $D_\sco$, as we now explain in more detail.
\subsection{Discretisation and method}

Several choices are available both for discretising the operators involved as well as 
minimising the criterion $\mathscr{J}$ defined by~\eqref{def:Lu}. 

\paragraph{Discretising the operator $T_a$.}
Let us quickly mention how we may synthetically produce data, \textit{i.e.}, how the operator $T_a$ is discretised.
The interval $[0,1]$ is evenly separated with $n$ points
$(t_i)_{0\leq i\leq n-1}$, $t_i =\frac{i}{n}$, 
with step $\Delta t = 1/n$.
A function $x$ is represented by the vector $X = (x(t_i))_{0 \leq i \leq n-1}$.
First,
for a constant kernel $k=1$ and any $a>0$, 
the discretised operator $\widetilde{T}_a$ may be computed as an approximation of $T_a$ defined by~\eqref{def:Abel}, thanks to the trapezoidal rule
\[
\begin{split}
 T_a x (t_i) = & \; \int_0^{t_i}(t_i -s)^{a-1}   x(s) \d s \\
 \left(\widetilde{T}_a   X \right)_i         = & \; \sum_{j=1}^i \frac{x(t_j)+x(t_{j-1})}{2} 
                      \int_{t_{j-1}}^{t_{j}}(t_i -s)^{a-1} \d s + O(\Delta t) \\
            = & \;  \sum_{j=1}^i \frac{x(t_j)+x(t_{j-1})}{2a} 
                      ((t_i -t_{j-1})^{a} - (t_i -t_{j})^{a}) + O(\Delta t)
            \; ,
\end{split}
\]
which leads to the corresponding matrix 
\[
(\widetilde{T}_a)_{i,j}
= 
\begin{cases}
\; \displaystyle \frac{(\Delta t)^a}{2a}
     \left( (i-j+1)^{a} - (i-j-1)^{a}  \right) & j<i \; ,\\
     \\
     \; \displaystyle \frac{(\Delta t)^a}{2a}  (i^a - (i-1)^a) & j=0, \; i \neq 0   \; ,\\
\\
\; \displaystyle \frac{(\Delta t)^a}{2a}   & j=i , \; i \neq 0 \; ,\\
\\
\; 0 & i=j=0, \text{ or } j > i.
\end{cases}
\]
This approximation is of order one as shown in \cite{diethelm1997} or \cite{li2011numerical}.
Thus, from the discrete operator $\widetilde{T}_a$, 
we compute $Y= \widetilde{T}_a X$,
to which we add a Gaussian normal noise of different standard deviations $\delta$ to obtain the data $Y^d$.

\paragraph{Discretising the derivative $D_\sco$.}
Instead of computing the finite difference of order $p$, we discretise the Hilbert scale.
For $\sco \in \N^*$, we define the matrix $B_\sco$ as the approximation of $(-\Delta)^\sco$ with the appropriate boundary conditions. In order to do so, we directly compute the matrix $B_\sco$ with the finite difference method, where the boundary conditions are enforced in the construction.

Several discretisation choices are possible depending on the wanted order. For instance, using discretisations that are all of order (at least) $O(\Delta t^2)$, here are the resulting matrices for $\sco = 1, 2, 3$.

For $\sco=1$, $B= - \Delta$, $x'(0)=0$ and $x(1)=0$,
\[
B_1 = (\Delta t)^{-2}
\begin{pmatrix}
2 & -2 &0 & \cdots &0 \\
-1& 2 &-1&\cdots &0\\
0& -1 &2 &-1 & 0\\
0 & 0 & \ddots &\ddots & \ddots \\
\end{pmatrix}
\; .
\]
For $\sco=2$, $B= \Delta^2$, $x''(0) = x'''(0) =0$ and $x(1)=x'(1)=0$,
\[
B_2 = (\Delta t)^{-4}
\begin{pmatrix}
2 & -4 &  2   & 0 & \cdots & \cdots &0 \\
 -2  &   5   & -4 & 1 &\cdots & \cdots & 0\\
1 & -4  &   6    & -4 & 1 & 0& 0 \\
0  & 1 & -4 & 6  & -4 & 1 &0 \\
0  &  0 & \ddots & \ddots & \ddots &\ddots& \ddots\\
 0 & 0& \cdots & 1 & -4 & 6 & -4\\
 0 & 0& \cdots & 0 & 1 & -4 & 7
\end{pmatrix}
\; .
\]
For $\sco=3$, $B= -\Delta^3$, with the corresponding boundary conditions,
\[
B_3 = (\Delta t)^{-6}
\begin{pmatrix}
2 & -6 &  6    & -2 & 0& 0& \cdots & \cdots &0 \\
-3 &   10   & -12 & 6 & -1 &0& \cdots & \cdots & 0\\
3 &   -12   & 19& -15 & 6&-1 &0 & \cdots & 0\\
-1 & 6  &   -15    & 20 & -15 & 6&-1 & 0 &0\\
0&-1 & 6  &   -15    & 20 & -15 & 6&-1 & 0 \\
\vdots & 0  &  0 & \ddots & \ddots & \ddots &\ddots& \ddots &0\\
\vdots &0 &   0   & -1& 6 & -15 & 6 &-15 & 20 \\
\vdots &0 &   0  & 0 & -1& 6 & -15 & 20 &-16 \\
\vdots & 0 &  0    & 0& 0&  -1& 6 & -14 & 26 \\

\end{pmatrix}
\; .
\]

\begin{remarque}
Note that the discretised form of $B_\sco$ for $\sco=1,2$ or $3$ no longer is symmetric, 
although the continuous operator is. This can certainly be circumvented by considering the weak formulation of the elliptic partial differential equation $B_\sco x = y$ and using (for instance) finite elements, but requires more involved computations which we believe make the numerical approach less straightforward, while not improving its efficiency. 

Note, however, that all we need for solving the problem is to compute fractional powers of $B_\sco^\ast B_\sco$ and not $B_\sco$ directly. Hence, the powers are also uniquely defined at the discrete level since the matrices we need to take powers of are all symmetric. 
\end{remarque}
To compute the fractional power of symmetric matrices, we use the Schur-Padé alogorithm developed in~\cite{higham2011}.

\paragraph{Minimisation of $\mathscr{J}$.}
Finally, we minimise the discrete function over $\R^n$
\[
\widetilde{\mathscr{J}} (X) = \|\widetilde{T}_a X -Y^d\|^2
                     + \alpha \|D_\sco^p X\|^2
                     \; ,
\]
where the norm is the Euclidean norm over $\R^n$. Its unique minimum $X_{\delta,\alpha,p}$ satisfies the so-called normal equations
\begin{equation}
   \label{eq:solmin}
    X_{\delta,\alpha,p} 
    = 
    \left(\widetilde{T}_a^T \widetilde{T}_a
    + \alpha \big(\widetilde{D}_\sco^T\widetilde{D}_\sco\big)^p\right)^{-1} \widetilde{T}_a^TY^d
    \; .
\end{equation}
For the numerical implementation, 
we choose $n$ small, 
i.e. $n =100$. 
For experimental signals, 
this value may reach several thousands or more, 
and in that case the minimum of a quadratic function 
can efficiently be obtained by the conjugate gradient method.
Both techniques have been implemented and give the same results.


\subsection{Recovering the theoretical rates}

Solving the inverse problem through the Tikhonov approach highly depends on the parameter $\alpha$. In this subsection, our purpose is to retrieve the theoretical orders of convergence given by Corollary~\ref{cor:order}. 

Hence, for illustration purposes, we here and only here choose the value of the regularisation coefficient $\alpha$ optimally, 
\textit{i.e.}, by estimating the best possible reconstruction error as follows:
\begin{algorithm}[H]
 \KwData{true signal $x$, noisy signal $y^\delta$}
 \KwResult{Solve the inverse problem  by minimising~\eqref{def:Lu}}
 \While{ $\alpha_m < \alpha < \alpha_M$ }{
  Solve \eqref{eq:solmin} for $\alpha$ \;
  Compute error = $\| x_{\delta,\alpha,p} -x \|$\;
  \If{error $<$ error $\text{opt}$}{
   Select $\alpha$ as $\alpha_{\text{opt}}$\;
   }
   }
\end{algorithm}
\noindent
Here, the regularisation parameter is searched for in $[\alpha_m, \alpha_M]$, where $\alpha_m$ is set at $10^{-16}$ and $\alpha_M$ chosen appropriately depending on the data (for a function of norm~$1$ we pick~$\alpha_M = 10^4$). 
The regularisation parameter space is explored incrementally with a logarithmic step, and is increased as long as the prediction no longer improves.
Of course, 
such a search cannot be implemented in a practical inverse problem, 
since this requires knowing the solution. 

We numerically illustrate the rate of convergence obtained for various values of $a$ and $p$. 
In particular, we highlight that the choice of the matrix $B_\sco$ is critical 
when the solution and its derivatives do not (properly) vanish at the boundary.

As an illustrating example, we reconstruct a Gaussian function centered at $x=\frac{1}{2}$ and with sufficiently small variance so that it numerically boils down to a compactly supported function. 
We also consider an off-center Gaussian function, for which the function (and its derivatives) do not vanish at $x=0$.
Those functions are both infinitely smooth.
However, only the centered Gaussian belongs to $\Xkp$ 
for any $\sco>0$ and any $p\geq 0$ (at the numerical level and for sufficiently small variance).

\paragraph{Rate of convergence.}

We solve the minimisation problem defined by~\eqref{def:Lu} 
for varying noise levels 
and show that convergence rates are close to the optimal ones as given by Corollary~\ref{cor:order}.

\begin{center}
\begin{figure}
\begin{adjustwidth}{-2cm}{2cm}  
    \begin{tabular}{ccc}
    \includegraphics[width=0.4\textwidth]{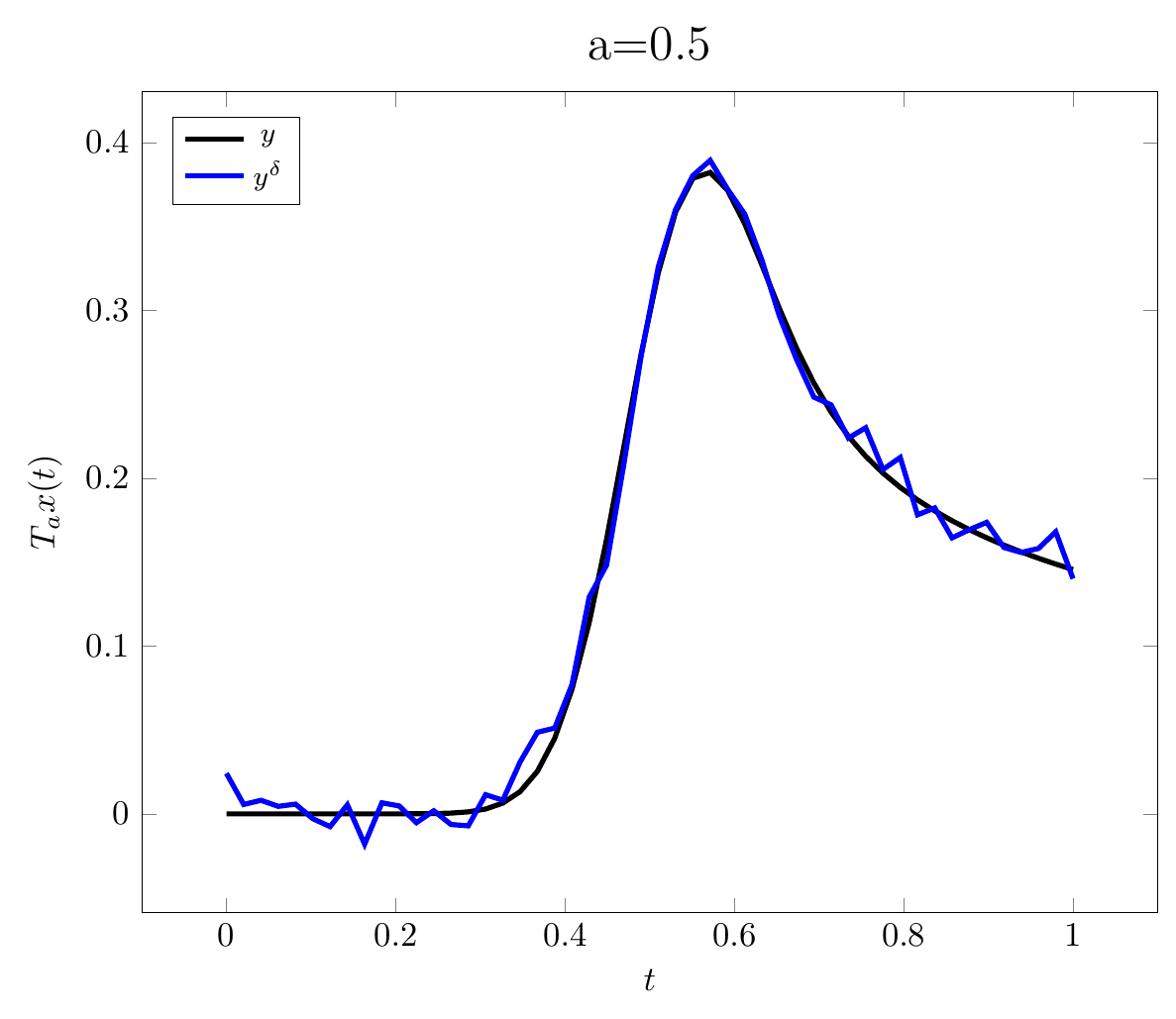}
    & \includegraphics[width=0.4\textwidth]{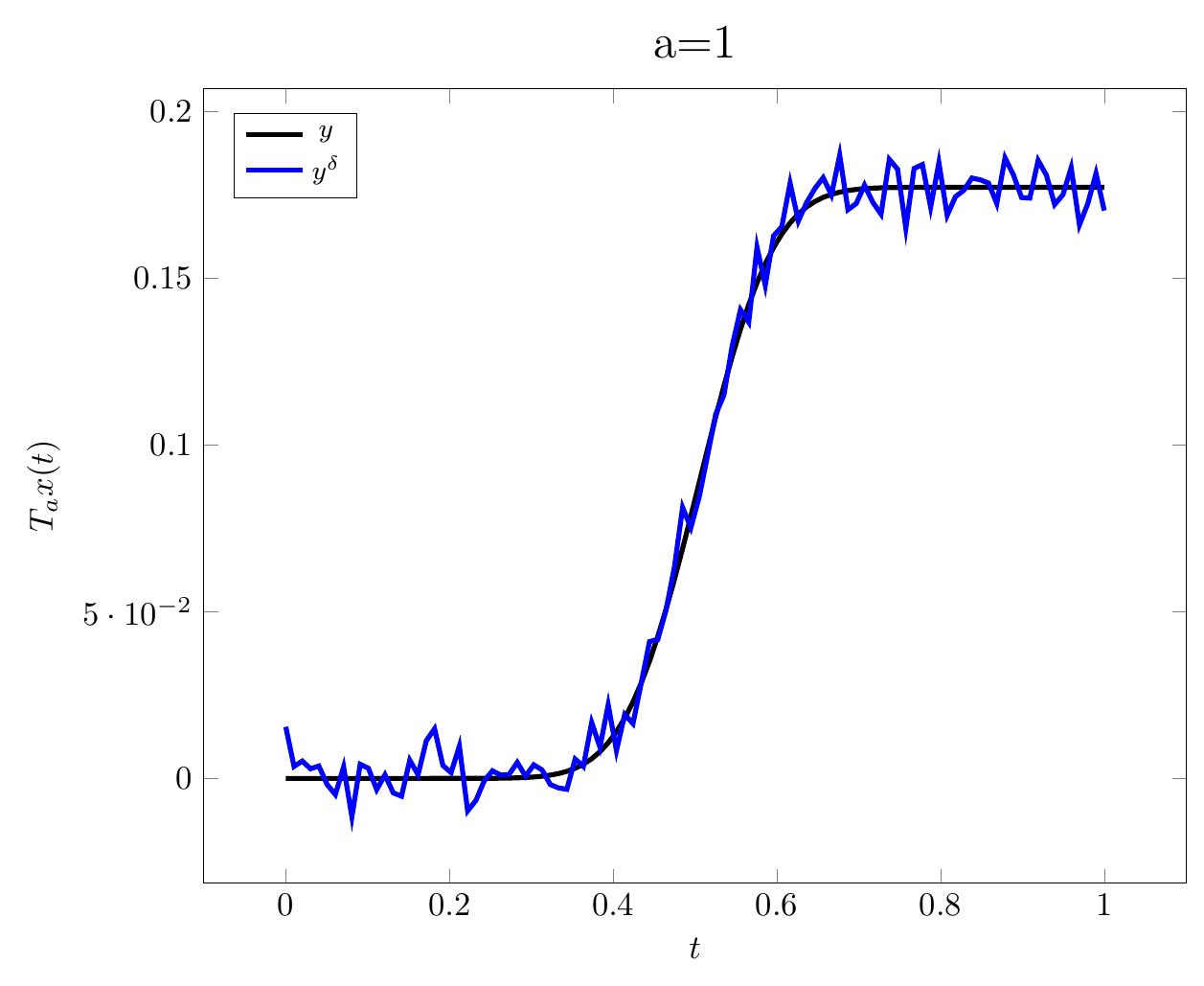}
    &\includegraphics[width=0.4\textwidth]{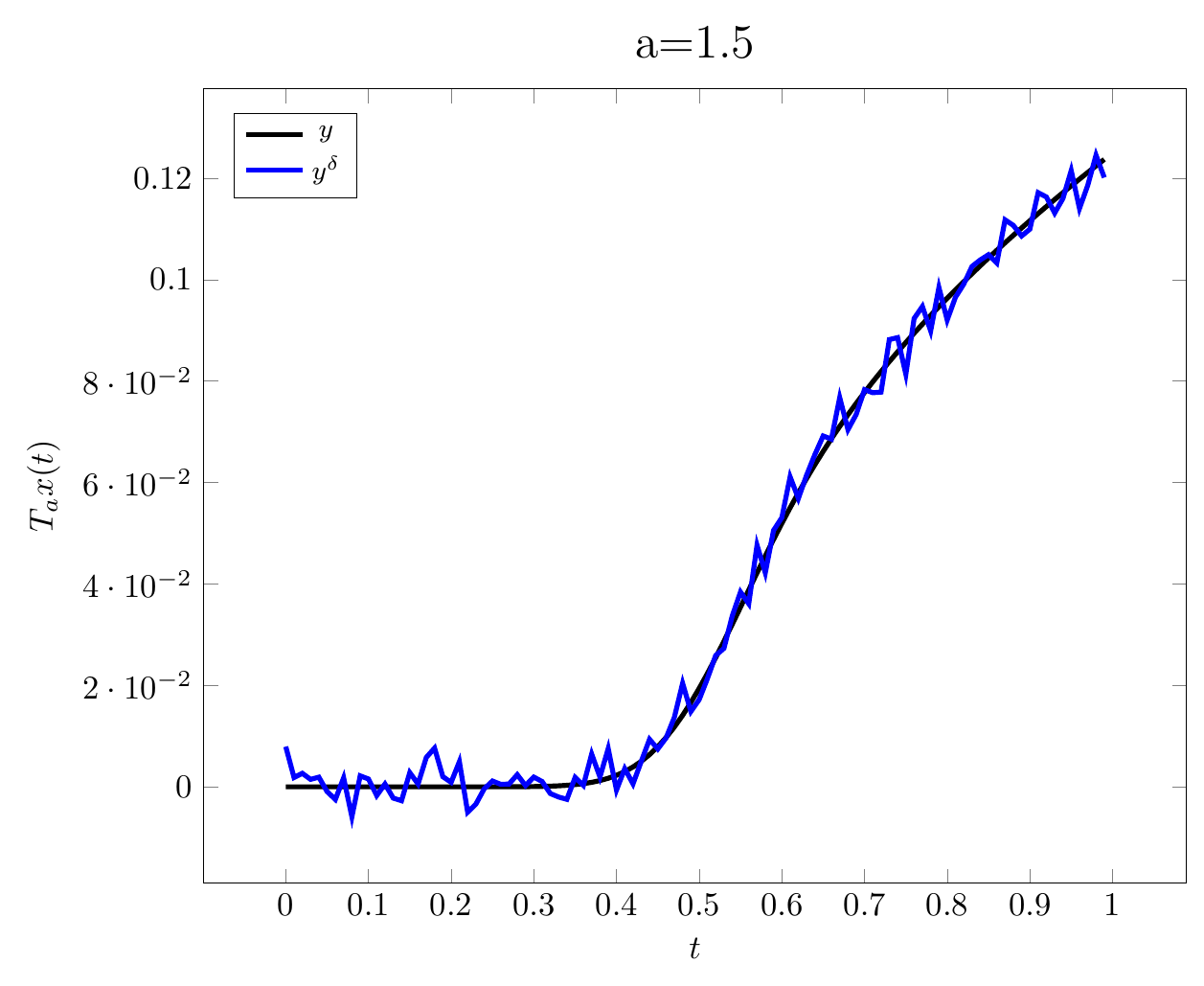} \\
    \includegraphics[width=0.4\textwidth]{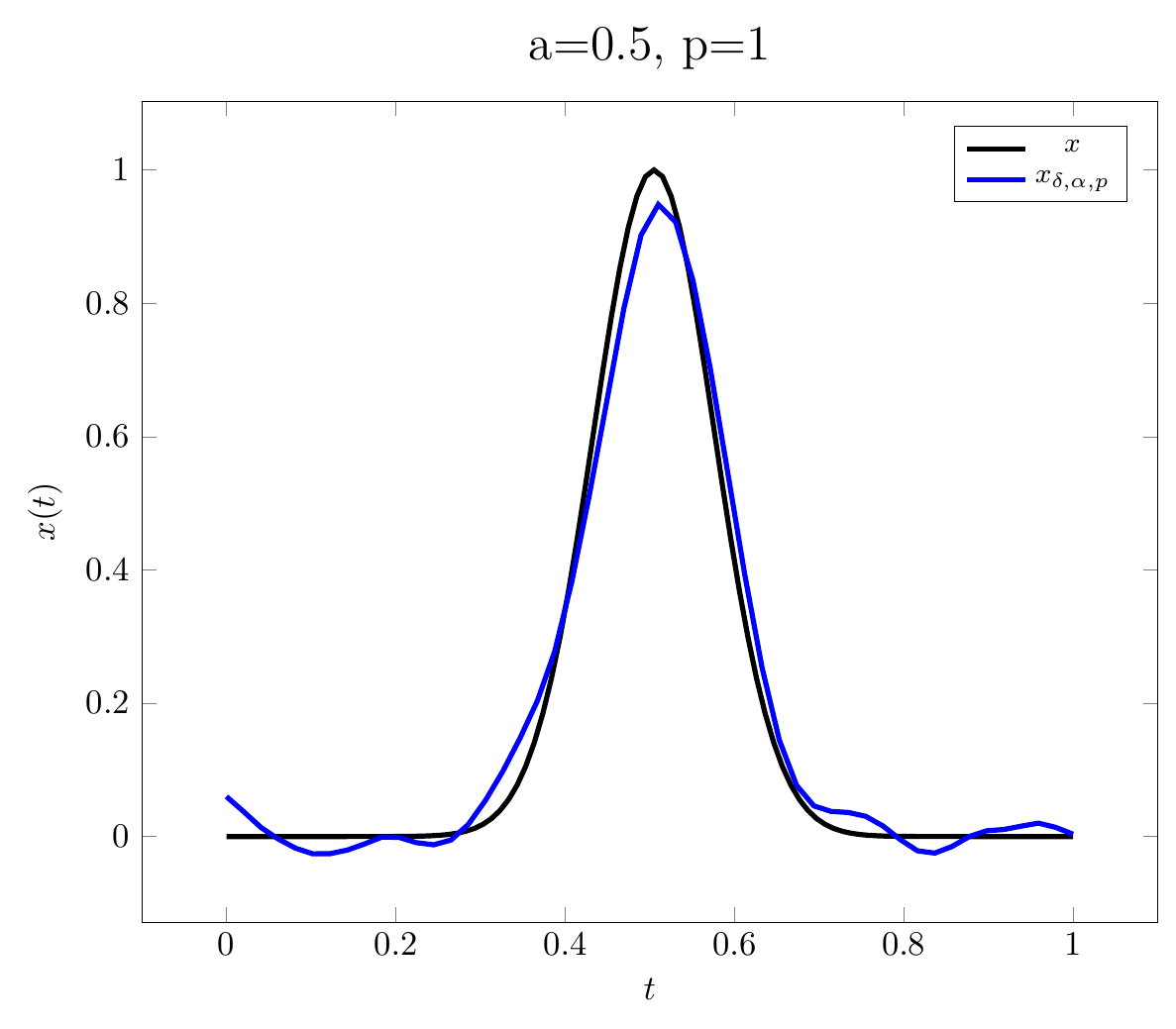}
    &\includegraphics[width=0.4\textwidth]{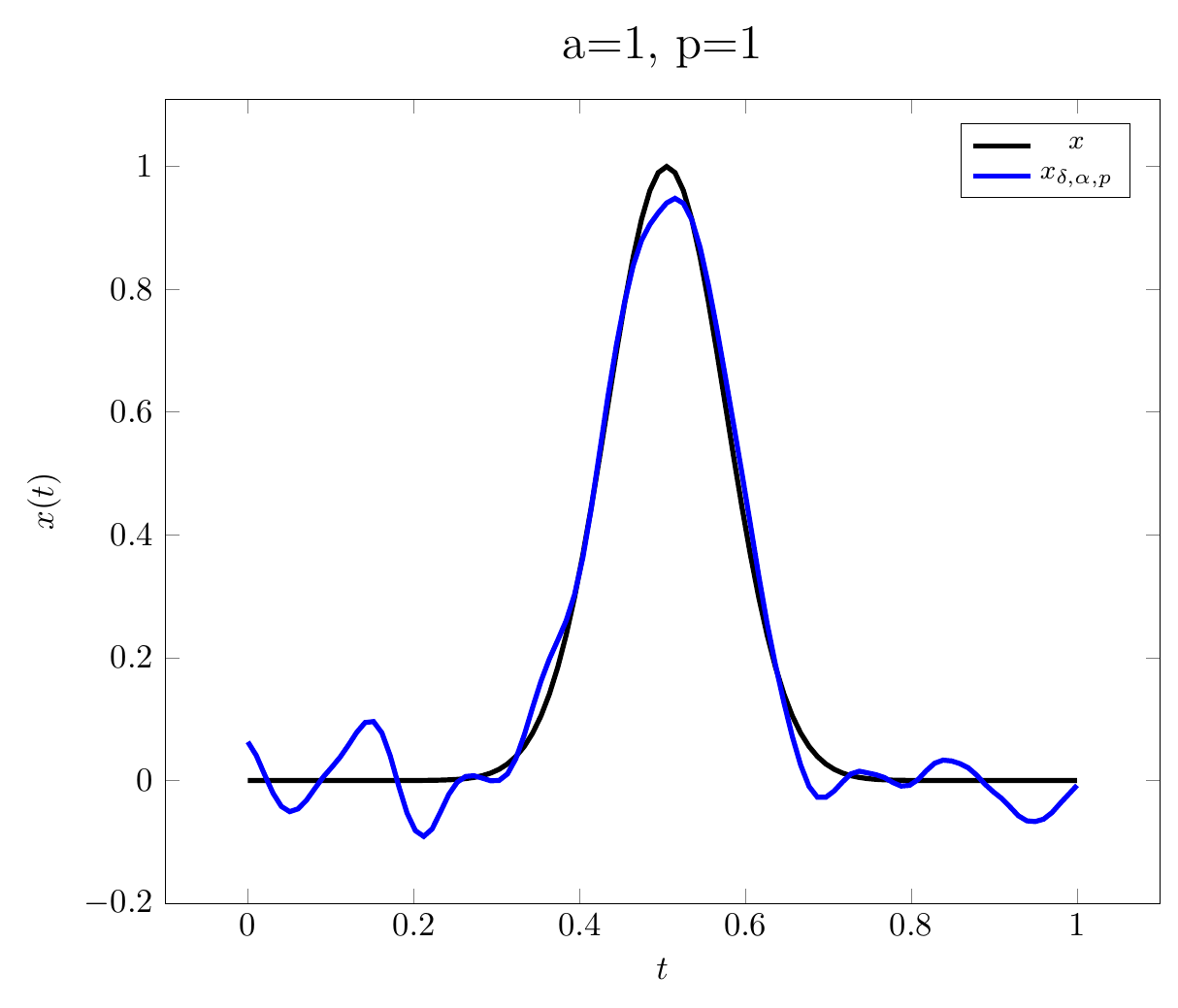}
    &\includegraphics[width=0.4\textwidth]{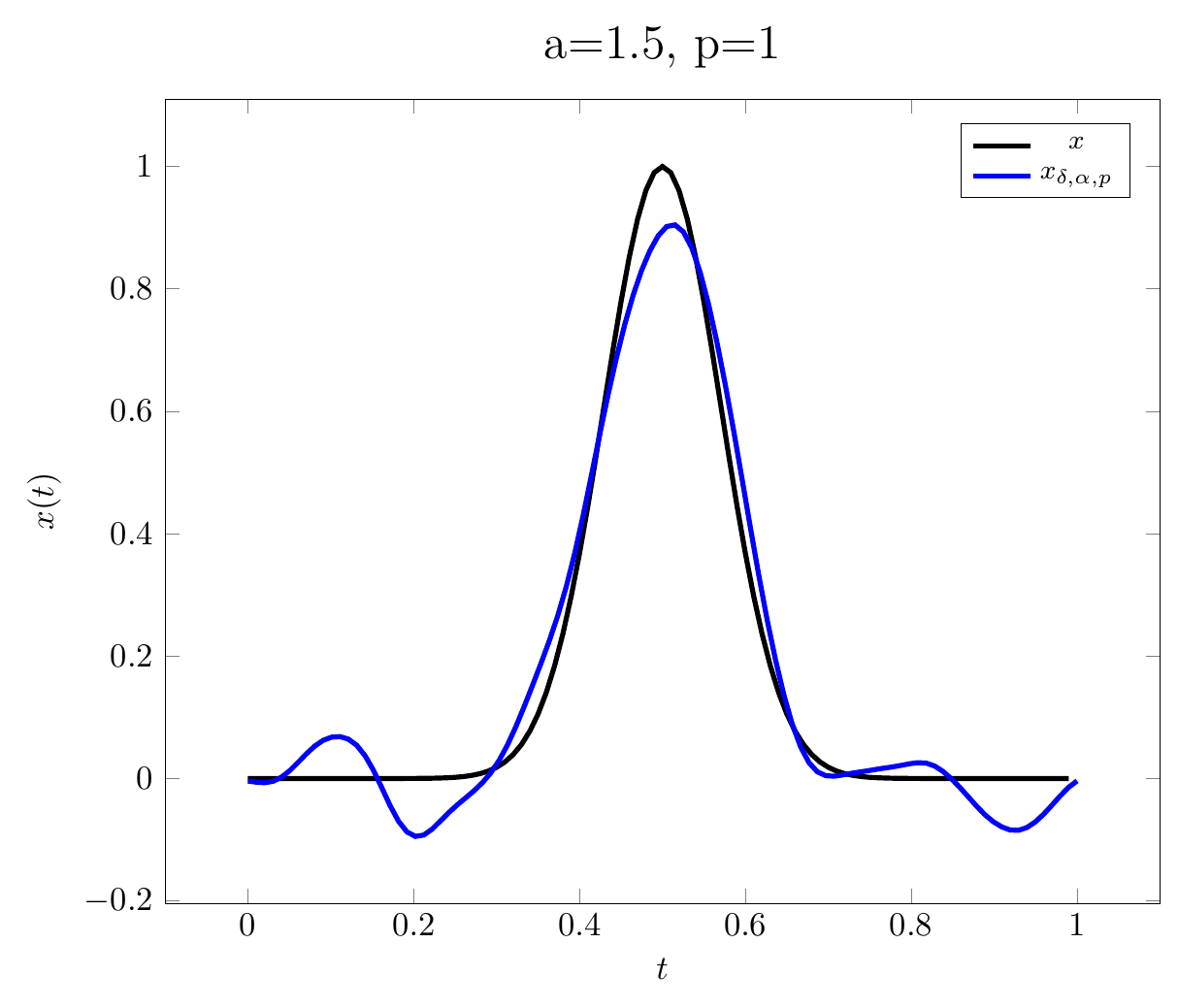}\\
    \includegraphics[width=0.4\textwidth]{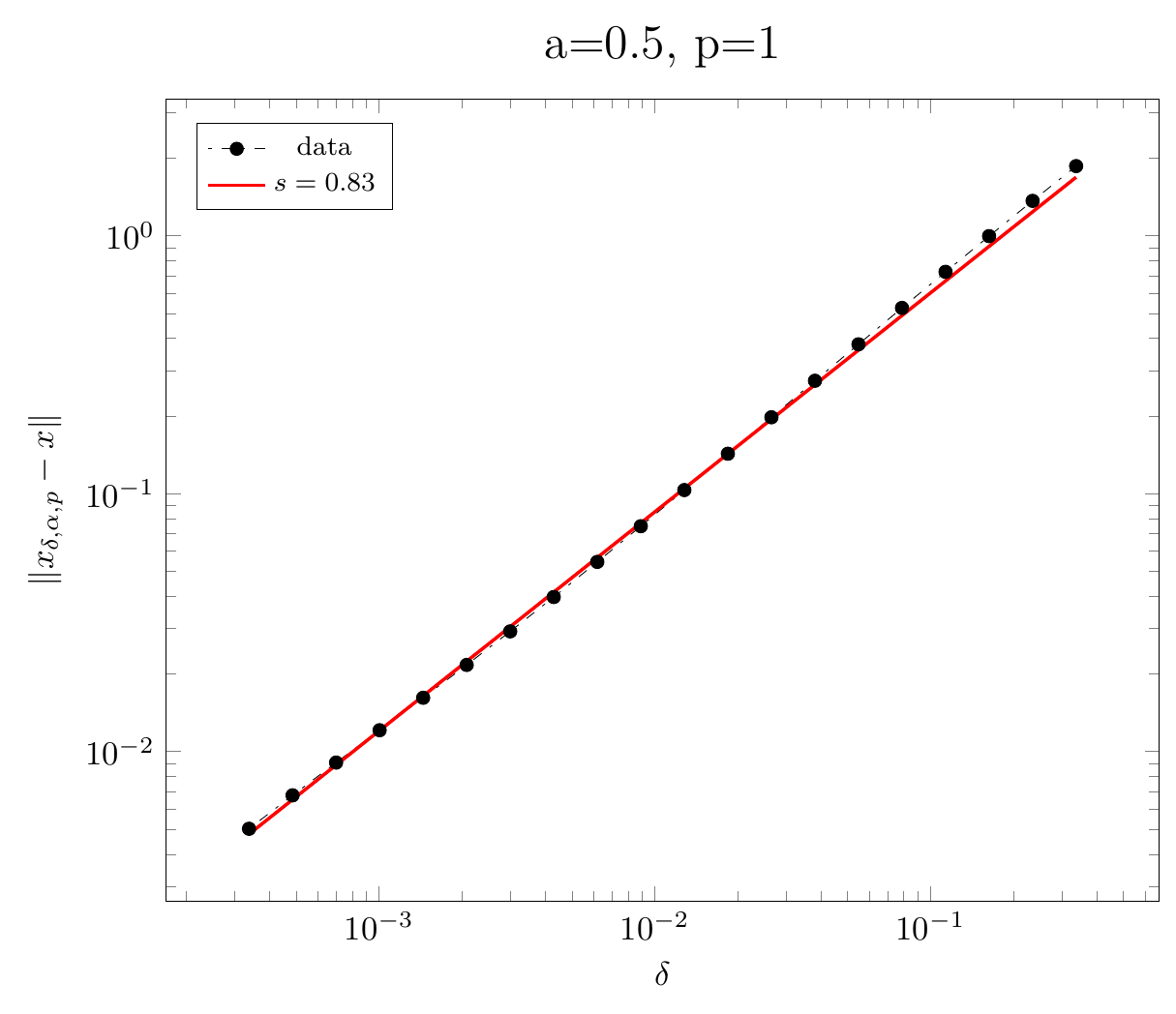}
    &\includegraphics[width=0.4\textwidth]{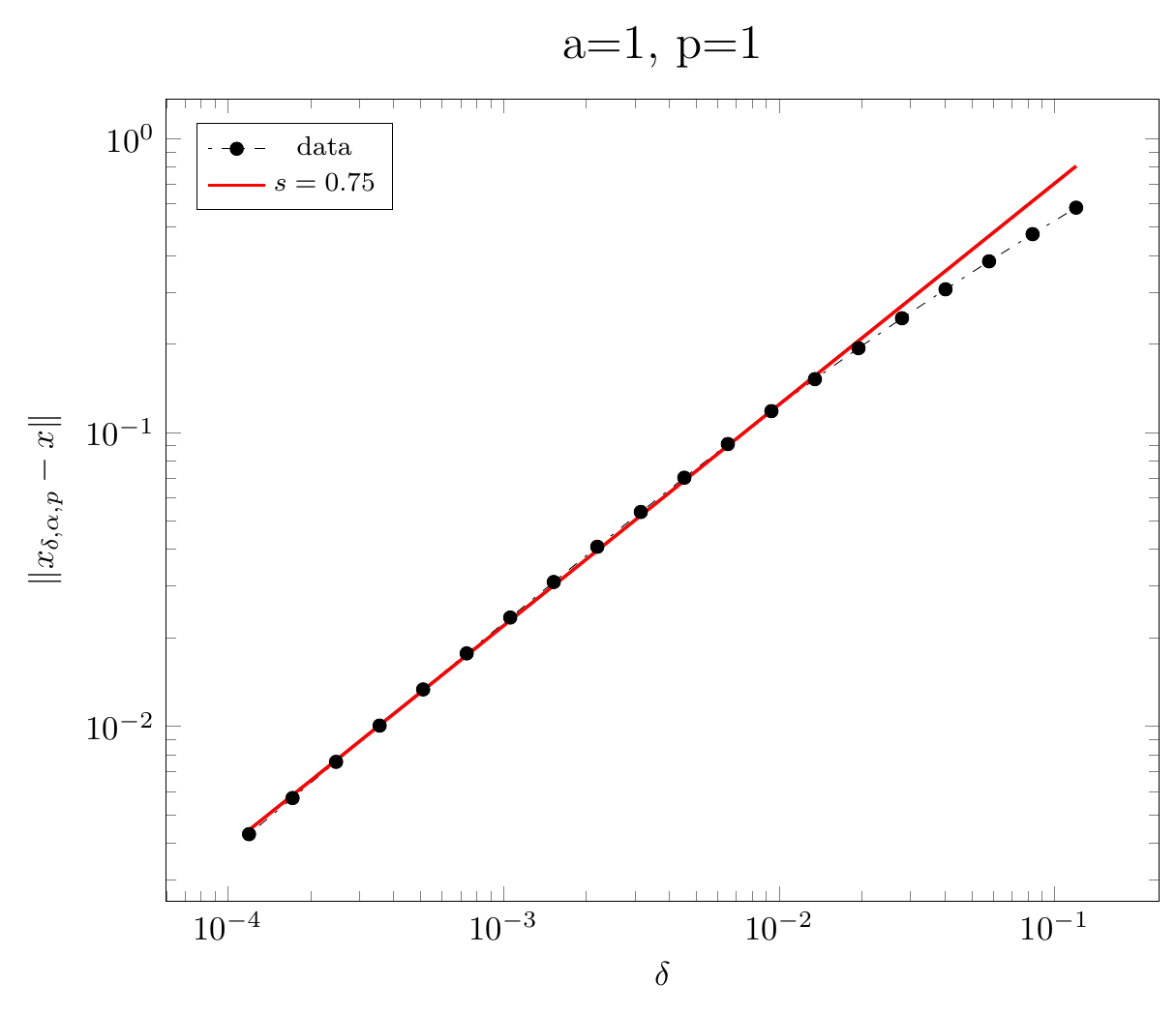}
    &\includegraphics[width=0.4\textwidth]{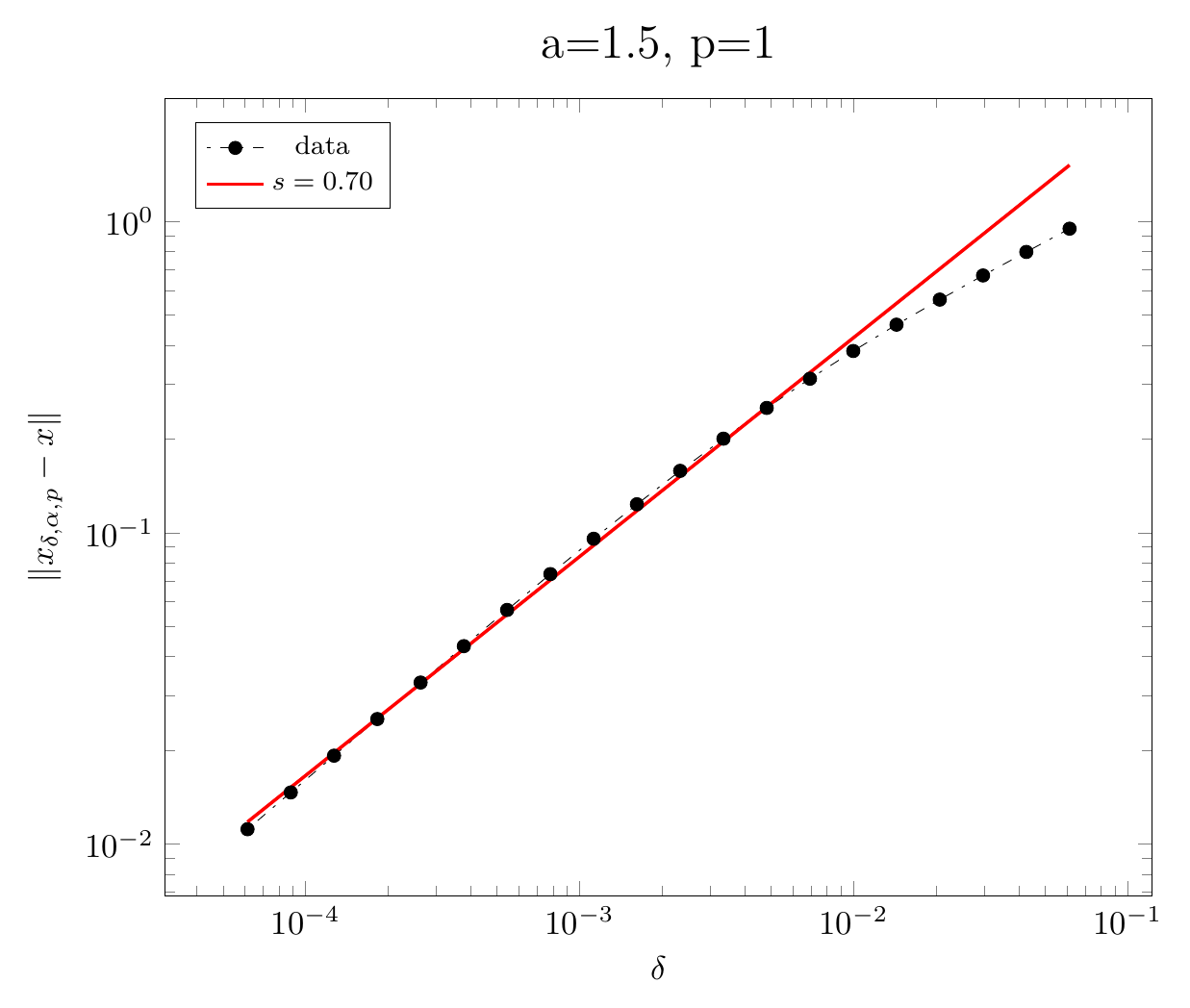}
    \end{tabular}
    \end{adjustwidth}
    \caption{
Examples of signals and their reconstructions for different values of~$a$. The first row shows the signal obtained by Abel transform of a Gaussian function for three different values of~$a$, in black the signals $y= T_ax$ without noise and in blue the noisy signals~$\yd$ with noise~$\delta = 0.05$.
The second row displays the reconstructions of the Gaussian from the noisy signals. 
Finally, in the third and last row, are plotted the reconstruction errors $\| x_{\delta,\alpha,p}- x\|$ as a function of the noise level $\delta$ (or equivalently the standard deviation of white Gaussian noise).
The smaller $a$, the better the reconstructed signal approximates the real solution~$x$ in $\log$-$\log$ scale.
The line of expected optimal slope $s$ is drawn in red. }
\label{fig:diffa}
\end{figure}
\end{center}

Figures~\ref{fig:diffa} shows the reconstruction of a Gaussian signal for different values of the order of ill-posedness~$a$ and of the order of penalisation $p$.
For a given level of noise, or standard deviation $\delta = 0.05$, 
the lower the order $a$, the better the reconstruction becomes.
Moreover, Figure~\ref{fig:diffa} shows that the optimal slope is attained, 
which for $p=1$ is $s=(2+a)/(2+2a)$, 
\textit{i.e.} $s \approx 0.833$, $s=0.75$ and $s=0.70$ for $a=0.5$, $a=1$ and $a=1.5$ respectively.

\paragraph{Example of saturation.}

We also aim at highlighting how the slope of convergence rates $s$ saturates. 
Indeed, 
assume that the unknown belongs to $X_{\sco,q}$ and let $p^* = \frac{q-a}{2}$.
As soon as $p$ is chosen higher than $p^*$, 
the convergence rate should no longer improve, as per Corollary~\ref{cor:order}. 

\begin{center}
\begin{figure}
\begin{adjustwidth}{-2cm}{2cm}  
    \begin{tabular}{ccc}
    \includegraphics[width=0.4\textwidth]{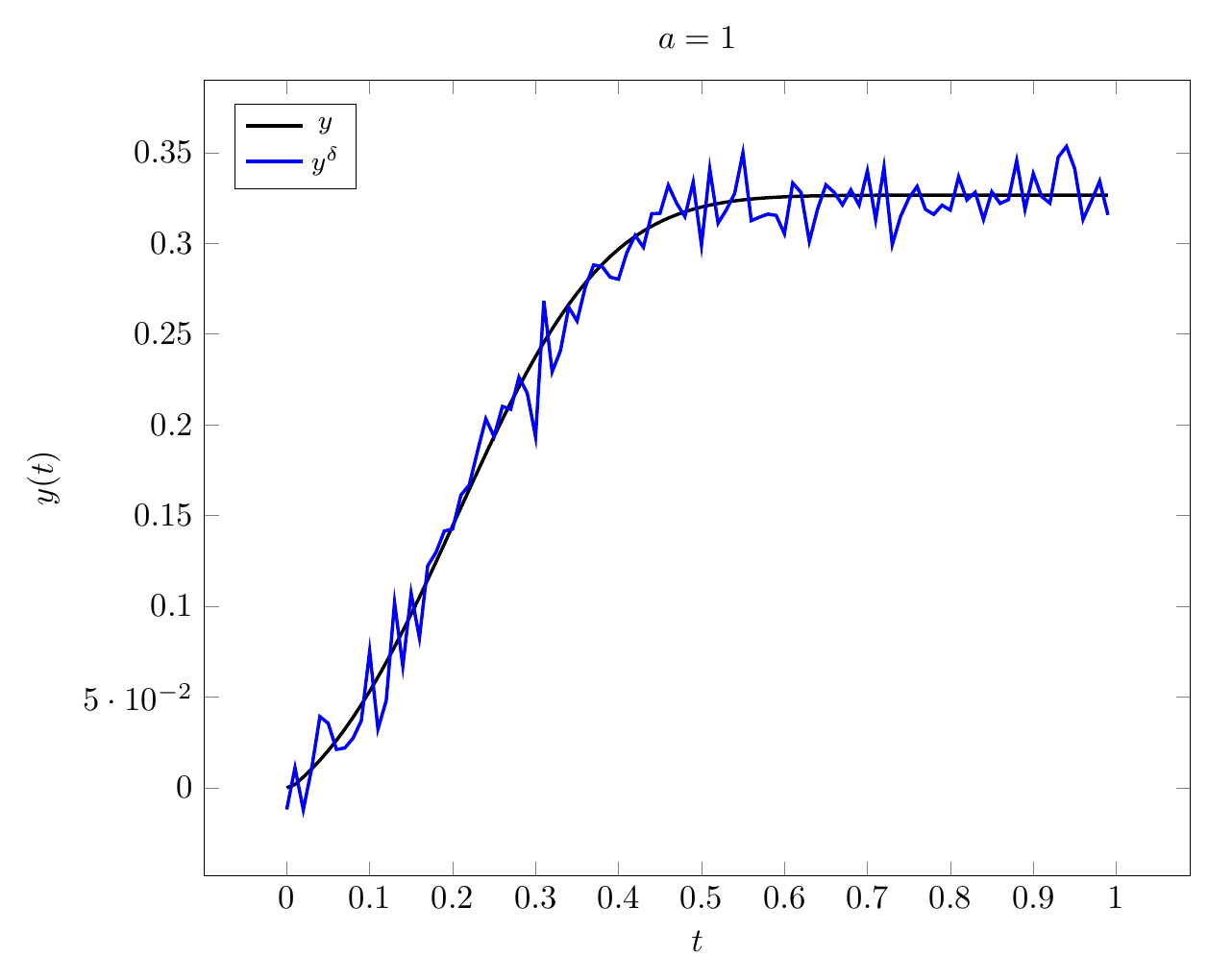}
    & \includegraphics[width=0.4\textwidth]{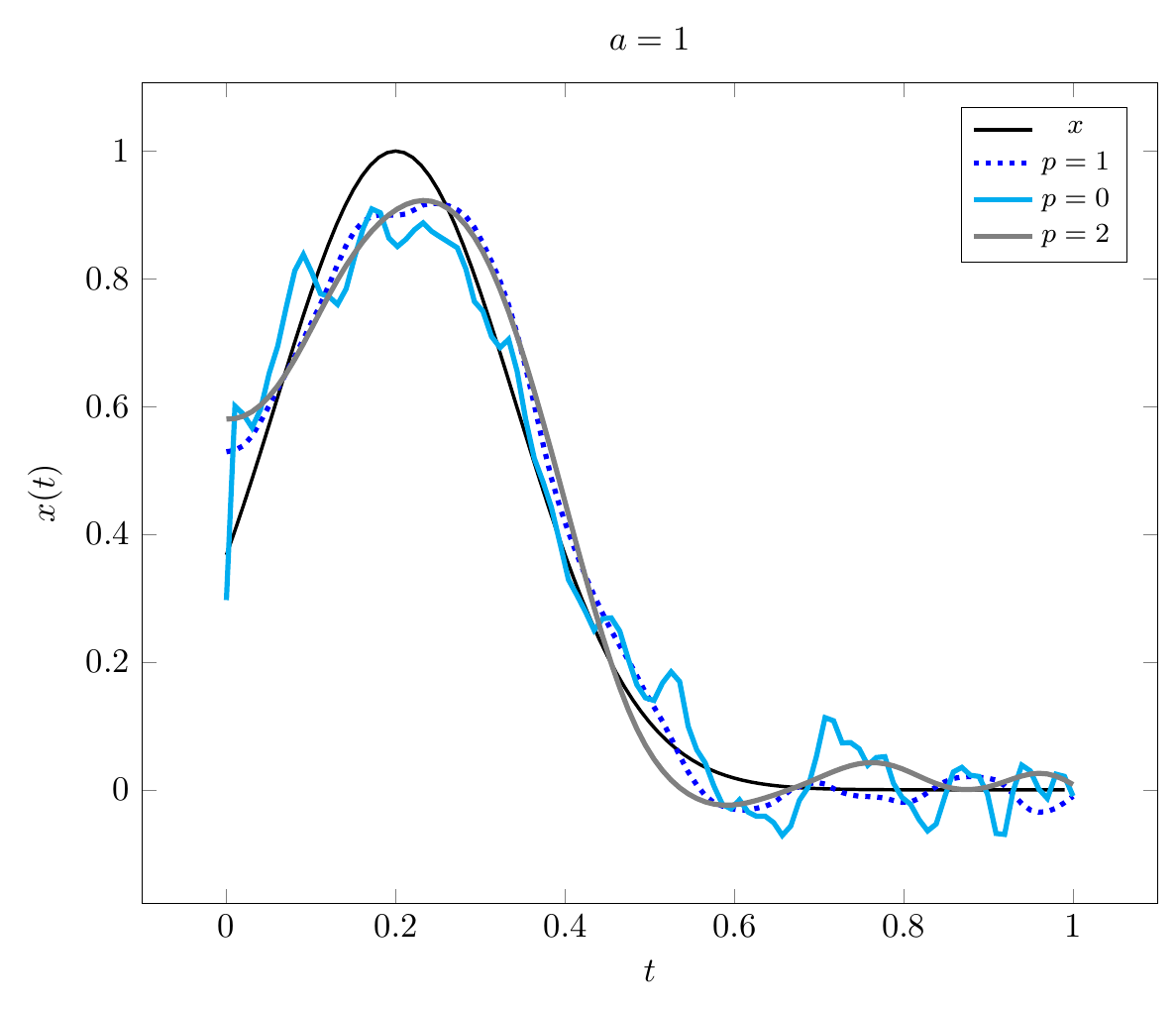}
    &\includegraphics[width=0.4\textwidth]{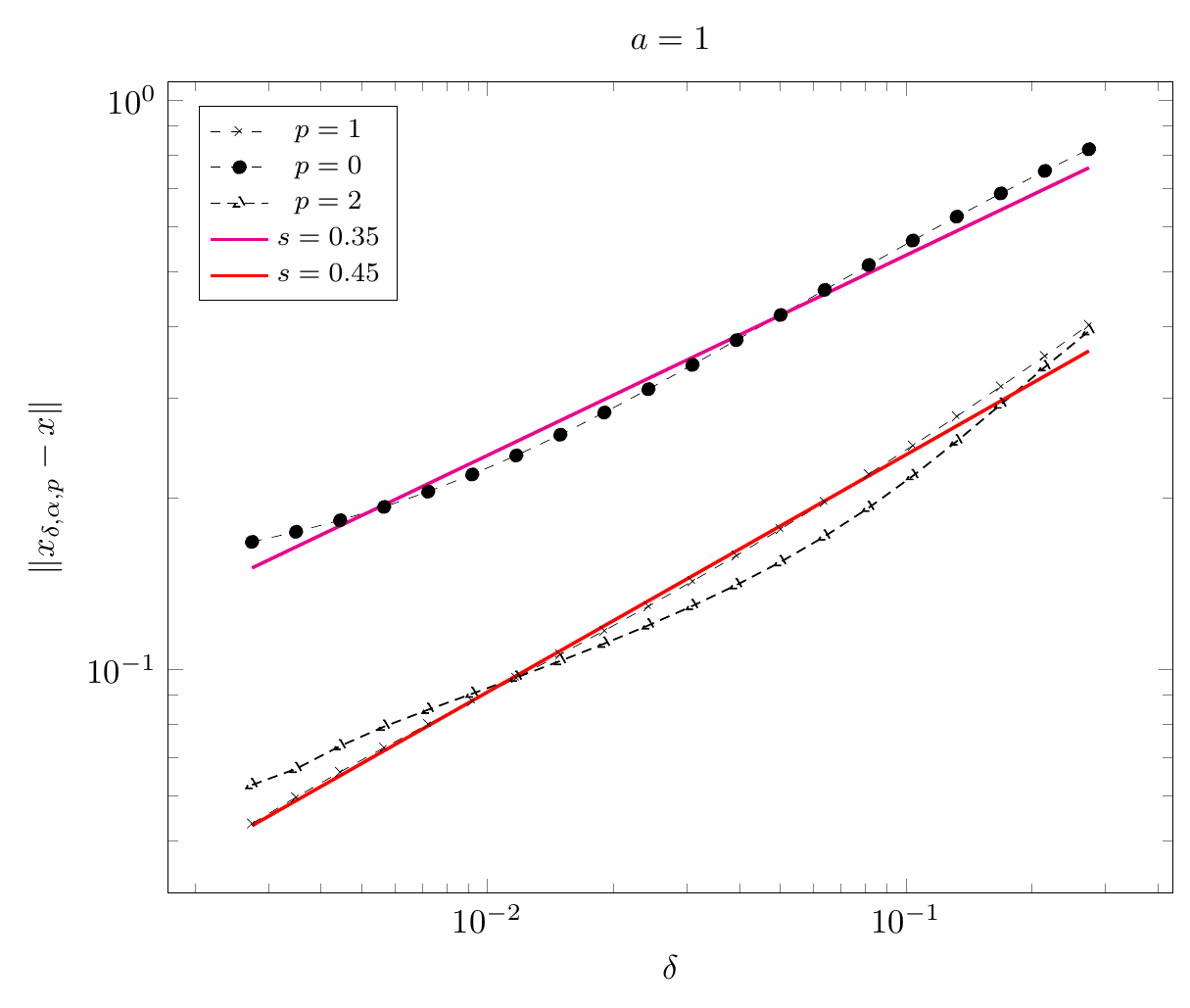} 
    \end{tabular}
    \end{adjustwidth}
    \caption{Example of slope saturation when the function $x$ 
    and its derivatives do not vanish at the boundary. 
    On the left is plotted the Abel transform of a Gaussian function centered on $t=0.2$.
    The noisy observation $\yd$ results from the addition of a white Gaussian noise with a standard deviation~$\delta =0.05$ to $y = T_ax$.
    In the middle is plotted the reconstrcution $x_{\delta,\alpha,p}$ for $p=0$, $p=1$ and $p=2$.
    The last graphic shows the error $\|x_{\delta,\alpha,p} -x \|$ as a function of $\delta$ in a $\log-\log$ scale.
    The slope $s$ for $p=1$ or $p=2$ is only $0.45$, which is close to the ideal $0.6$,
    instead of the slope of $0.75$ obtained in Figure~\ref{fig:diffa}.}
    \label{fig:saturation}
\end{figure}
\end{center}

Thus, for $a=1$, 
if we choose as function $x$ to be an off-center Gaussian, 
then $x(1)=0$ but $x'(0) \neq 0$, 
and this function therefore belongs to $X_{1,q}$ for any $q<3/2$. In particular it is in $X_{1,1}$.
However, it does not belong to $X_{1,q}$ for $q>3/2$.
We deduce that $p^*=(3/2-1)/2 = 1/4$ 
and for all $p \geq 1/4$ the slope remains at $q/(a+q) = 0.6$. 
This saturating phenomenon for the order of convergence is illustrated by Figure~\ref{fig:saturation}.
We notice that the slope is slightly below than the expected $0.6$.

\begin{remarque}
According to Remark~\ref{rem:commute}, 
when the operators commute, the saturation of the convergence rate $s$ is obtained with a smaller $p^*$, \textit{i.e.}, $p^*= q/2-a$.
In our example, the continuous operators $(S_a^*S_a)^{1/2}$ and $D_r$ commute.
On the other hand, the discrete operators lose this property.
\end{remarque}

\paragraph{The importance of the chosen derivative operator.}
As already mentioned, if the function $x$ is not compactly supported, the right choice of the derivative operator becomes crucial.
To illustrate this phenomenon, 
we again pick the off-center Gaussian. 

For $a=1.5$ (hence $\sco=2$), 
such a function belongs to $X_{2,q}$ for all $q<5/2$, and does not belong to $X_{2, q}$ for $q>5/2$
since $x(1)=x'(1)=0$, but $x''(0) \neq 0$.
Then, $p^* =1/2 (5/2-1.5) = 1/2$ and for $p \geq p^*$ the slope of convergence is $q/(a+q) = 0.62$.

\begin{center}
\begin{figure}
\begin{adjustwidth}{-2cm}{2cm}  
    \begin{tabular}{ccc}
    \includegraphics[width=0.41\textwidth]{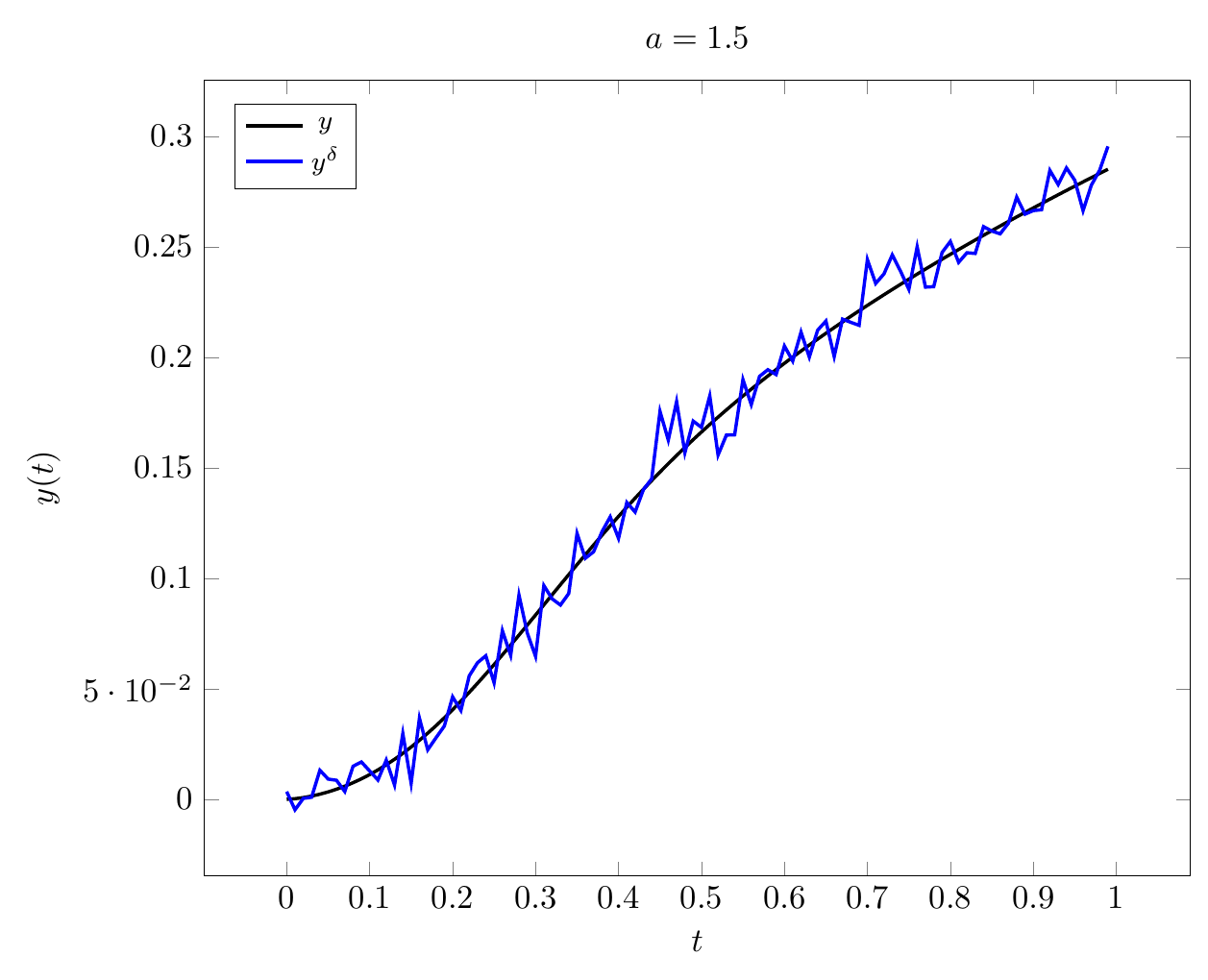}
    & \includegraphics[width=0.4\textwidth]{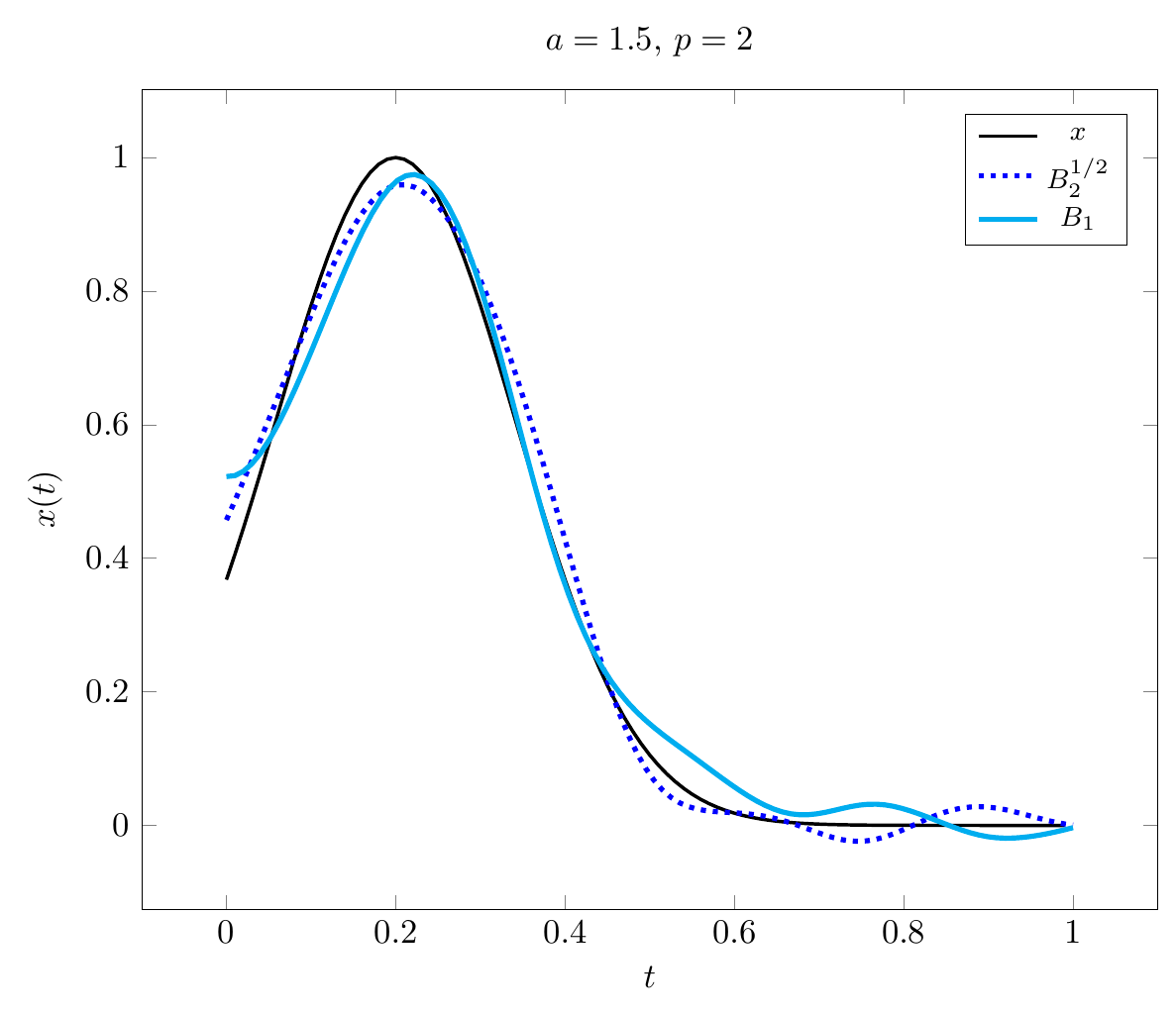}
    &\includegraphics[width=0.4\textwidth]{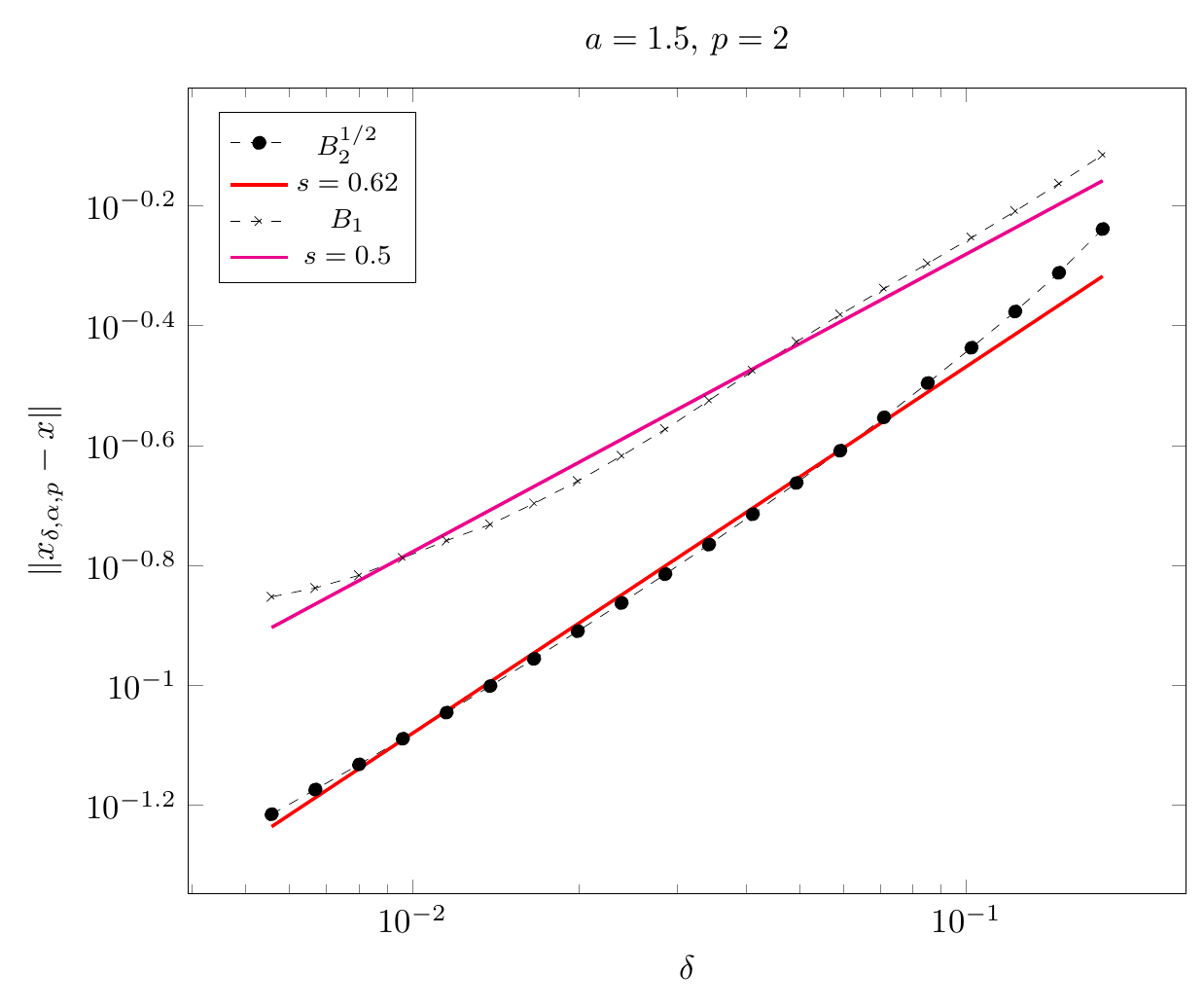} 
    \end{tabular}
    \end{adjustwidth}
    \caption{Effect of the chosen matrix $D_\sco$ and $B_\sco$. 
    For $a=1.5$ and $p=2$, and a standard deviation $\delta =0.05$,
    we compare the reconstructions obtained with either the matrix $B_2^{1/2}$ or $B_1$.
    On the left is plotted the Abel transform for $a=1.5$ of a Gaussian function centered on $0.2$. 
    The noisy observation $\yd$ results from the addition of a white Gaussian noise with a standard deviation~$\delta =0.05$ to $y = T_ax$.
    In the middle is plotted the reconstruction $x_{\delta,\alpha,p}$ for $p=2$ and two different derivative operators.
    On  the right is shown the error $\|x_{\delta,\alpha,p}-x\|$ as a function of the noise $\delta$ in a $\log-\log$ scale, and their slope $s$.
    The optimal rate of convergence $s=0.62$ is only obtained with $B_2^{1/2}$,
    which is the matrix associated to the right Hilbert scale for the inverse problem.}
    \label{fig:operator}
\end{figure}
\end{center}

We choose to compare for $p=2$ the effect of choosing either the matrix $B_1$ or $B_2^{1/2}$.
In Figure~\ref{fig:operator}, we observe that the matrix $B_1$ leads to a reconstruction $x_{\delta,\alpha,p}$ which must satisfy the condition $x'(0)=0$,
whereas such condition is not enforced with the choice of $B_2^{1/2}$.
We also notice that the slope of convergence is optimal for $B_2^{1/2}$, 
but saturates at $0.5$ for $B_1$.
This confirms that the matrix $B_2^{1/2}$ offers a better reconstruction,
and in that case the values taken by the unknown at the boundaries play an important role. We hence numerically confirm that regularising with a high order can be less effective if the operator is not chosen appropriately.

Moreover, we note that the solution becomes increasingly sensitive to the a priori regularity parameter $p$ as it becomes large. Hence, even if formally taking $p$ large cannot be harmful, it leads to numerical instabilities. This is an incentive to choose $p$ as optimally as possible depending on the problem under study.

\subsection{Example in stereology}

We now propose to apply our method to an example in stereology, with the aim of proceeding as in a real experimental situation. The model proposed in~\cite{evans2017} or~\cite{jakeman1975abel} reads
\[
y(t) = \sqrt{t} \; \int_0^t \frac{x(s)}{(t^2-s^2)^{1/2}}  \d s
\; ,
\]
which can be rewritten in the form
\begin{equation}
\label{ker:stereo}
y(t) = \int_0^t k(t,s) (t-s)^{a-1} x(s) \d s, \qquad k(t,s)=\frac{\sqrt{t}}{\sqrt{t+s}}
\; .
\end{equation}
Note that this kernel still falls short of satisfying our regularity hypotheses (as it already did with those of~\cite{gorenflo1991abel}). 
Indeed, it can be checked that $k$ does not satisfy the condition of~\eqref{cond:triangle} since $g(t,s) =  \; (k(t,t) -k(t,s) ) (t-s)^{a-1}$ is not even in $H^1(0,t)$.
We shall see that the method nonetheless works efficiently.

Upon using the trapezoidal rule, the discretisation of $T_a$ reads
\[
(\widetilde{T}_a)_{i,j}
= 
\begin{cases}
\; \displaystyle \frac{(\Delta t)^a}{2a} \frac{\sqrt{i}}{\sqrt{i+j}}
     \left( (i-j+1)^{a} - (i-j-1)^{a}  \right) & j<i \; ,\\
     \\
     \; \displaystyle \frac{(\Delta t)^a}{2a}  (i^a - (i-1)^a) & j=0, \; i \neq 0   \; ,\\
\\
\; \displaystyle \frac{(\Delta t)^a}{2a} \frac{1}{\sqrt{2}}   & j=i , \; i \neq 0 \; ,\\
\\
\; 0 & i=j=0 \text{ or } j>i \; .
\end{cases}
\]
We consider an initial vector $X$ of very large size, 
much larger than the reconstruction sample, namely $N \gg n$. 
We then compute $Y = (y_i)_{0 \leq i \leq N-1}$, 
to which a white noise of unknown standard deviation is added, 
chosen in the interval $[0.01,0.1]$. 
After sub-sampling the signal, 
we obtain $Y^d~=~(y^\delta_i)_{0 \leq i \leq n-1}$,
from which we reconstruct the signal $x_{\delta,\alpha,p}$.

In order to solve this inverse problem, we pick the smoothing operator associated to $a=0.5$, 
which are respectively the square root of $B_1$ for $p=1$ and $B_1$ for $p=2$. 
Since we do not have access to the true data $x$ or the noise level $\delta$,
we follow the discrepancy principle as an a posteriori rule  to select the parameter $\alpha$~\cite{engl1996regul}.
More precisely, we first assume that the signal 
$(y_i)_{0 \leq i\leq i \, \text{max}}$ 
is null up to some known time $t<t_{i_{\text{max}}}$. Then, there is only noise and the noise level $\delta$ may therefore be estimated as the average of $(y_i)_{0 \leq i\leq i \, \text{max}}$.
Then, the regularisation parameter $\alpha$ is chosen 
so that the error is in the same range 
as the expected noise level $\delta$.

\begin{figure}
    \centering
    \includegraphics[scale=0.4]{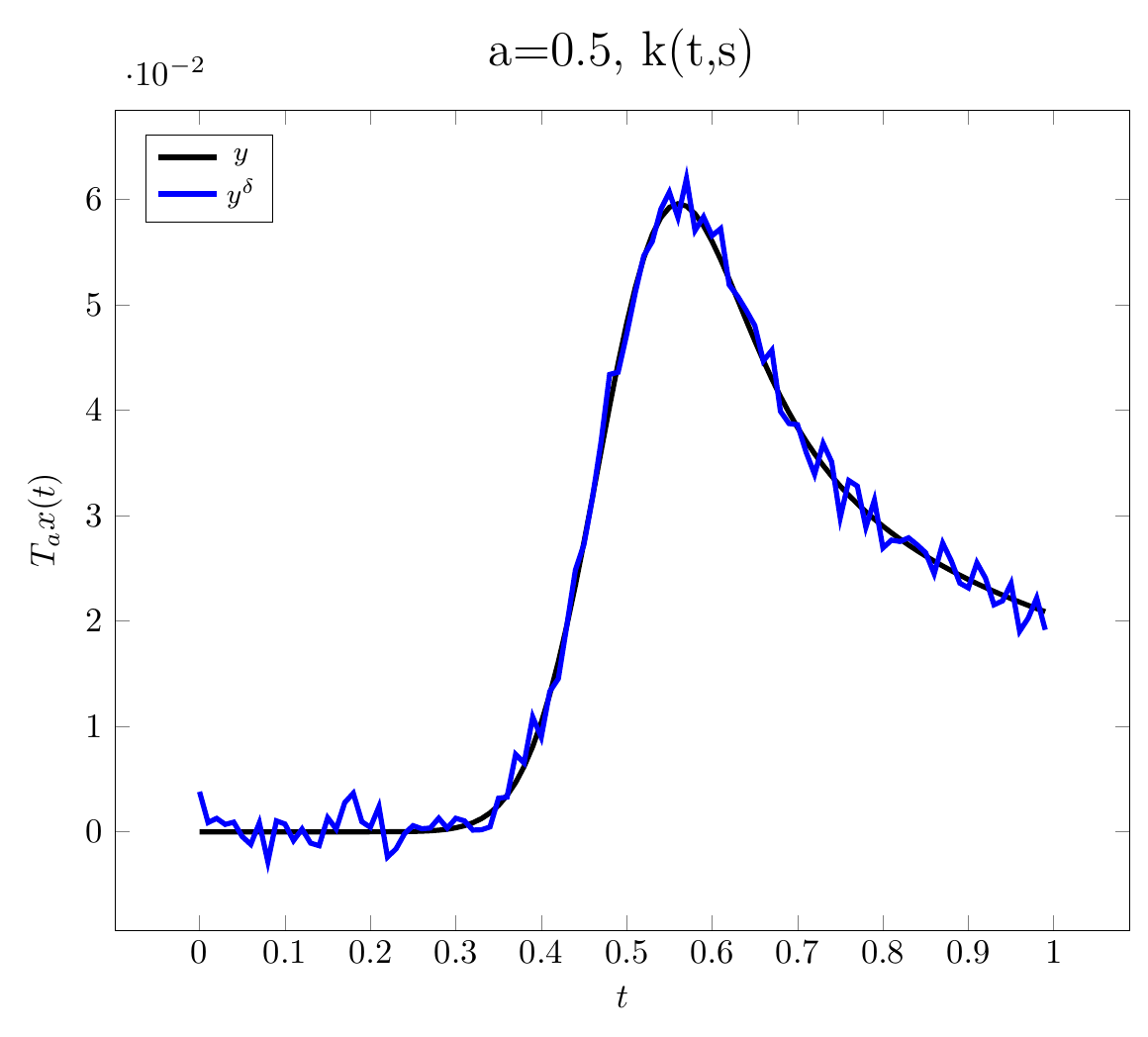}
    \includegraphics[scale=0.4]{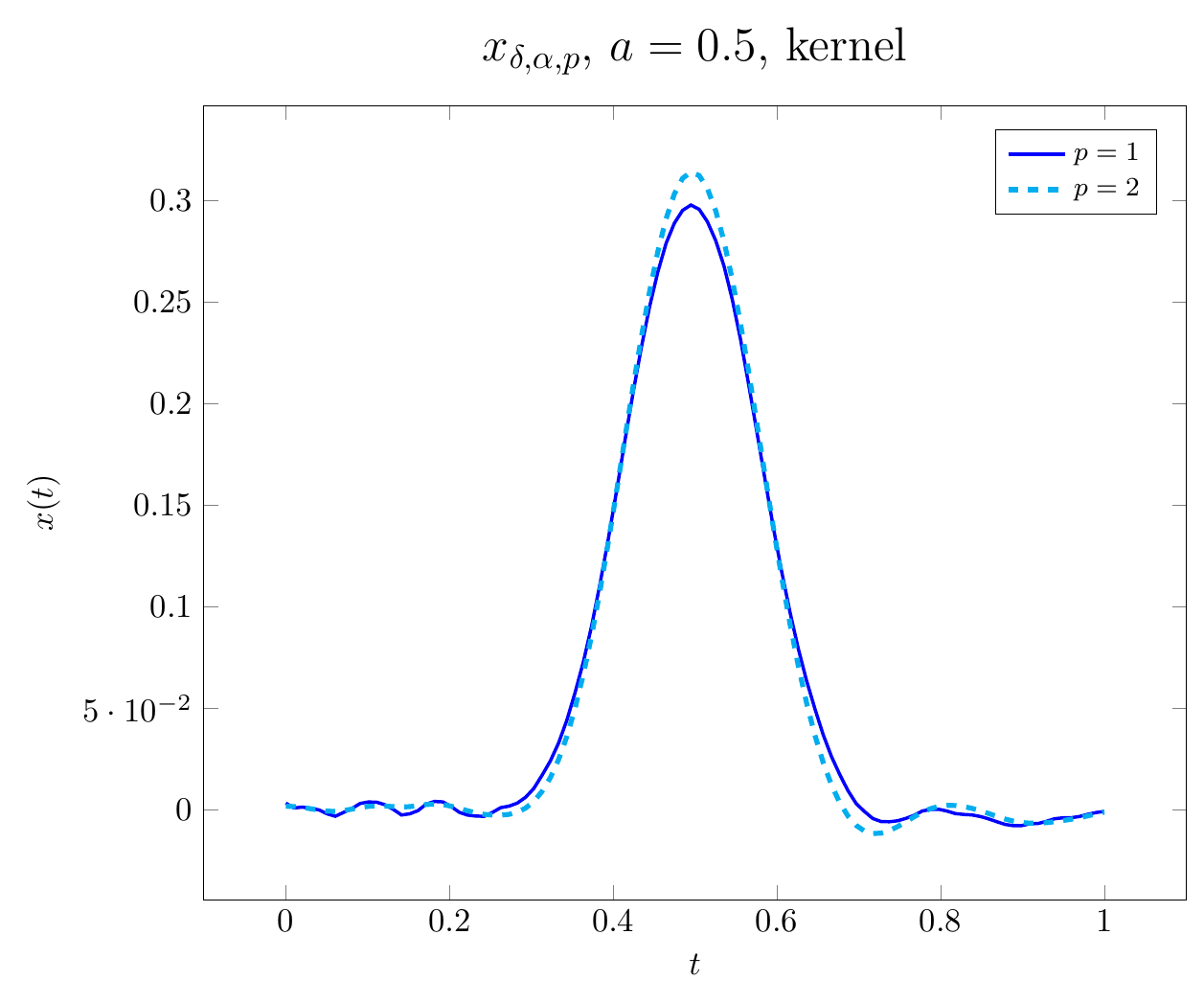}
    \caption{Reconstruction by Tikhonov regularisation and a posteriori rule, for $a=0.5$, kernel $k$ given by~\eqref{ker:stereo}. On the left we plot the Abel transform of a Gaussian function centered on $t=0.5$. The measurement $y =T_a x$ is corrupted with an additive white Gaussian noise with a standard deviation~$\delta =0.05$.
    On the right we show the reconstruction $x_{\delta,\alpha,p}$ obtained with $p=1$. }
    \label{fig:a0.5kernel}
\end{figure}

\begin{algorithm}[H]
 \KwData{noisy signal $y^\delta$}
 \KwResult{Solve the inverse problem  by minimisation of~\eqref{def:Lu}}
 Compute $\delta$ as the average of $(y^\delta_i)_{0\leq i \leq i \text{max}}$ \;
 \While{ $\alpha_m < \alpha < \alpha_M$ }{
  Solve \eqref{eq:solmin} for $\alpha$ \;
  Reconstruct $T x_{\delta,\alpha,p}$ \;
  Compute error = $\| T x_{\delta,\alpha,p} - y^\delta \| $ \;
  \If{error $\sim $ $\delta$}{
   Select $\alpha$ as $\alpha_{\text{opt}}$\;
   }
   }
\end{algorithm}
Figure~\ref{fig:a0.5kernel} shows the reconstruction $x_{\delta,a,p}$
for $p=1$ or $2$. 
As expected, the reconstruction is smoother for $p=2$.
Even if the parameter $\alpha$ is not optimal, 
the reconstruction method for $x_{\delta,\alpha,p}$ remains efficient when combined to a posteriori rules dealing with the unknown level of noise. 

\paragraph{Acknowledgments.}
The authors are grateful to Nikolaos Roidos for the insightful exchanges about the Heinz-Kato inequality and his work~\cite{roidos2019heinz}.



\bibliographystyle{unsrt}
\bibliography{biblioRegul}

\appendices
\section{The analytic case}
\label{annexe:proof}

We here elaborate on making condition~\eqref{cond:Retafini} more explicit at the expense of requiring more smoothness for $k$.
Assuming that for each $t \in (0,1)$, $k(t,\cdot)$ is analytic around~$t$ with radius of convergence at least~$t$,
we may write for all $(s,t) \in \Omega$,
\[
k(t,s) = \sum_{n \in \N} \frac{(-1)^n a_n(t)}{n!} (t-s)^n
\; .
\]
We introduce a family of operators indexed by $n \in \N^*$ by letting $a_n(t):= \partial_s k(t,t)$, and
\[
A_{n, \poc} x(t) = a_n(t) \int_0^t (t-s)^{n-\poc-1} x(s) \d s
\; .
\]
We give here an other computation of the splitting of $T_a$ between the main term comparable to $S_a$ and a compact perturbation of this term.
For $x \in L^2$,
\[
\begin{split}
    T_R x (t)= & \;  \sum_{n=1}^{+\infty} \frac{(-1)^n a_n(t)}{n!} S_{a+n} x(t) \\
             = & \sum_{n=1}^{+\infty} \frac{(-1)^n a_n(t)}{n!}
    \frac{\Gamma(\sco+n) \Gamma(\res)}{\Gamma(\sco) \Gamma(\res + \poc) \Gamma(n) \Gamma(\poc) }
    S_{n-\poc} S_{a+\poc} x(t)
    \; .
\end{split}
\]
Then, noticing that $a_n(t) S_{n-\poc} = A_{n, \poc}$, we have
\[ T_R x(t) = \sum_{n=1}^{+\infty} \frac{(-1)^n}{n!}
    \frac{\Gamma(\sco+n) \Gamma(\res)}{\Gamma(\sco) \Gamma(\res + \poc) \Gamma(n) \Gamma(\poc) }
    A_{n,\poc} S_{a+\poc} x(t)
    = R_{a,\poc} S_{a+\poc} x(t)
    \; .
\]
Moreover, Stirling's approximation yields as $n \rightarrow +\infty$.
\[
\begin{split}
 \frac{\Gamma(\sco+n) \Gamma(\res)}{\Gamma(\sco) \Gamma(\res + \poc) \Gamma(n) \Gamma(\poc) } \sim 
\frac{n^\sco}{\Gamma(\sco) }
\frac{\Gamma(\res)}{\Gamma(\res + \poc)\Gamma(\poc)}  \sim C n^{\sco}
\; .   
\end{split}
\]
In fact, the expression of $R_{a,\poc}$ in the form~\eqref{def:Raeta} is equivalent to the one above, as can be seen from an explicit calculation.
For  $0<\poc< \sco-a$ small enough, a sufficient condition for the condition $\| R_{a,\poc}\| < +\infty$ to hold then is
\begin{equation}
    \sum_{n=1}^{+\infty} \frac{n^r}{n!} \|A_{n,\poc}\| < +\infty 
    \; . \label{cond:k}
\end{equation}
Note that the analyticity of $k$ means 
condition~\eqref{cond:k} implicitly assumes that the operators $A_{n, \poc}$ are well-defined and bounded as operators from $L^2$ onto $L^2$, for $\poc$ small enough.
Condition~\eqref{cond:k} is still formulated in a general and abstract way,
but can be easily checked in practice. 
Let us make it more explicit in the following cases:
\begin{itemize}
    \item Let us suppose that the functions $a_n$ are bounded and $\|a_n\|_\infty = O(n^\gamma) $,
    then assuming that $A_n$ is an Hilbert-Schmidt operator,
    we may estimate its norm  for $a \notin \N$:
   \[ 
   \begin{split}
       \|A_n\|^2 \leq & \; \|A_n\|_{HS}^2\\
               = &\; \int_0^1 a_n(t)^2 \int_0^t (t-s)^{2n-2\poc-2} \d s \d t \\
               \leq & \; \frac{1}{(2n-2\poc-1)(2n-2\poc)} \; \|a_n\|_\infty^2\\
               = & \; O( n^{2\gamma-2} ) \; .
   \end{split}
   \]
   Those integrals are well defined as long as we choose $\poc<1/2$.
   The series $ \sum \frac{n^{r+\gamma -1}}{ n!}$ converges, 
   which ensures that condition~\eqref{cond:k} is met.
   \item Let us suppose that $a_n(t) = b_n t^{-\beta}$,
   for $\beta < 1$, and $a\notin \N$,
   \[ 
   \begin{split}
       \|A_n\|^2 \leq &\; b_n^2 \int_0^1 t^{-2\beta}\int_0^t (t-s)^{2n-2\poc-2}  \d s \d t \\
               \leq & \; \frac{b_n^2 }{(2n-2\poc-1)} \;  \int_0^1 t^{-2\beta + 2n-2\poc-1} \d t \; ,
   \end{split}
   \]
   These integrals are finite if $\poc<1$ is taken sufficiently small so that $\beta + \poc< 1$.
  If the series $\sum b_n \frac{n^{r-1}}{n!}$ converges, 
   condition~\eqref{cond:k} is met.
   %
 
   
\end{itemize}

\end{document}

%% file: command.tex
\let\d\undefined

\newcommand{\d}{\; \mathrm{d}}

\newcommand{\1}{{\mathchoice {\rm 1\mskip-4mu l} {\rm 1\mskip-4mu l}{\rm 1\mskip-4.5mu l} {\rm 1\mskip-5mu l}}}

\newcommand{\R}{\mathbb{R}}

\newcommand{\N}{\mathbb{N}}
\newcommand{\C}{\mathbb{C}}

\newcommand{\X}{X} 
\newcommand{\Y}{Y} 

\newcommand{\yd}{y^\delta} 
\newcommand{\xda}{x^{\delta,\alpha}}

\newcommand{\sco}{r} 
\newcommand{\res}{\omega} 
\newcommand{\Xkp}{X_{\sco,p}}

\newcommand{\tXkp}{\widetilde{X}_{\sco,p}}

\newcommand{\tXkpa}{\widetilde{X}_{\sco,p-2\sco}}

\newcommand{\Ta}{T_a}

\newcommand{\poc}{\gamma}